\documentclass[12pt,a4]{amsart}\usepackage{amssymb,euscript,array,bbm}
\newtheorem{Lemma}{Lemma}[section]\newcommand{\bel}{\begin{Lemma}}\newcommand{\eel}{\end{Lemma}}
\newtheorem{Proposition}[Lemma]{Proposition}\newcommand{\bprop}{\begin{Proposition}}\newcommand{\eprop}{\end{Proposition}}
\newtheorem{Theorem}[Lemma]{Theorem}\newcommand{\bthe}{\begin{Theorem}}\newcommand{\ethe}{\end{Theorem}}
\def\bull{\vrule height .9ex width .9ex depth -.1ex }
\newcommand{\bpr}{~\\{\em Proof.~}}\def\epr{$\bull$\\}
\newtheorem{Remark}[Lemma]{Remark}\newcommand{\beR}{\begin{Remark}\rm}\newcommand{\eeR}{\end{Remark}}
\newtheorem{Definition}[Lemma]{Definition}\newcommand{\bed}{\begin{Definition}}\newcommand{\eed}{\end{Definition}}
\newtheorem{Example}[Lemma]{Example}\newcommand{\bex}{\begin{Example}\rm}\newcommand{\eex}{\end{Example}}
\newtheorem{Corollary}[Lemma]{Corollary}\newcommand{\bcor}{\begin{Corollary}}\newcommand{\ecor}{\end{Corollary}}
\newtheorem{DefProp}[Lemma]{Definition-Proposition}\newcommand{\bdp}{\begin{DefProp}}\newcommand{\edp}{\end{DefProp}}
\newtheorem{Property}[Lemma]{Property}\newcommand{\bpro}{\begin{Property}}\newcommand{\epro}{\end{Property}}

\newcommand{\beq}{\begin{equation}}\newcommand{\eeq}{\end{equation}}
\newcommand{\bem}{\begin{displaymath}}\newcommand{\eem}{\end{displaymath}}
\newcommand{\beqa}{\begin{eqnarray}}\newcommand{\eeqa}{\end{eqnarray}}
\newcommand{\bee}{\begin{enumerate}}\newcommand{\eee}{\end{enumerate}}
\newcommand{\bei}{\begin{itemize}}\newcommand{\eei}{\end{itemize}}
\newcommand{\bet}{\begin{tabular}{cccccccc}}\newcommand{\eet}{\end{tabular}}
\newcommand{\bpm}{\begin{pmatrix}}\newcommand{\epm}{\end{pmatrix}}
\newcommand{\bM}{\begin{matrix}}\newcommand{\eM}{\end{matrix}}
\newcommand{\ber}{\begin{array}{l}}\newcommand{\eer}{\end{array}}

\def\bl{\langle}\def\br{\rangle}
\def\isom{\xrightarrow{\sim}}
\def\into{\stackrel{i}{\hookrightarrow}}\def\norm{\stackrel{\nu}{\to}}

\newcommand{\quotient}[2]{{\left.\raisebox{0.4ex}{$#1$}\!\!\middle/\!\!\raisebox{-0.4ex}{$#2$}\right.}}

\newcommand{\li}{~\\ $\bullet$ }
\newcommand{\ls}{~\\ $\star$ }

\def\cL{\mathcal{L}}\def\cM{\mathcal{M}}\def\cN{\mathcal{N}}\def\cO{\mathcal{O}}

\def\cm{{\frak m}}

\def\one{{1\hspace{-0.1cm}\rm I}}
\def\zero{\mathbb{O}}

\def\C{\mathbb{C}}
\def\k{\mathbbm{k}}\def\mN{\mathbb{N}}
\def\P{\mathbb{P}}
\def\R{\mathbb{R}}

\def\Z{\mathbb{Z}}

\def\al{\alpha}\def\be{\beta}
\def\ep{\epsilon}

\def\tA{\tilde{A}}\def\tB{\tilde{B}}\def\tC{{\tilde{C}}}

  \def\tX{{\tilde{X}}}

\def\empty{\varnothing}
\def\suml{\sum\limits}\def\oplusl{\mathop\oplus\limits}
\def\cupl{\mathop\cup\limits}

\def\cMv{{\cM^{\vee}}}\def\cNv{{\cN^{\vee}}}

\def\smin{\setminus}\def\sset{\subset}

\def\dr{determinantal representation}

\def\mg{maximally generated}\def\Mg{Maximally generated }
\def\CCS{{$C'/C$-saturated }}\def\XXS{$X'/X$-saturated }
\def\CM{Cohen-Macaulay }\def\cOn{\cO_{(\k^n,0)}}

\newcommand{\bin}[2]{\binom{#1}{#2}}

\title[D\MakeLowercase{eterminantal representations}]{D\MakeLowercase{eterminantal representations of
singular hypersurfaces in} $\P^n$}
\author[D.K\MakeLowercase{erner}]{D\MakeLowercase{mitry} K\MakeLowercase{erner}}
\address{Department of Mathematics, University of Toronto, 40 St. George Street, Toronto, Canada}
\email{dmitry.kerner@gmail.com}
\date{\today}
\author[V.V\MakeLowercase{innikov}]{V\MakeLowercase{ictor} V\MakeLowercase{innikov}}
\address{Department of Mathematics, Ben Gurion University of the Negev, P.O.B. 653, Be'er Sheva 84105, Israel.}
\email{vinnikov@math.bgu.ac.il}
\thanks{Part of the work was done during postdoctoral stay of D.K. in Mathematics Department of Ben Gurion
University, Israel.
 Both authors were supported by the Israel Science Foundation.
\\The authors thank A.Beauville, R.O.Buchweitz, G.M.Greuel and E.Shustin for numerous important discussions.
}
\subjclass[2000]{Primary 14M12;  Secondary 14H50; 14H99;}

\keywords{determinantal hypersurfaces, arithmetically Cohen-Macaulay sheaves, hyperbolic polynomials}

\begin{document}
\maketitle  \setcounter{secnumdepth}{6} \setcounter{tocdepth}{1}
\begin{abstract}
A (global) determinantal representation of projective hypersurface $X\sset\P^n$
is a matrix whose entries are linear forms in
homogeneous coordinates and whose determinant defines the hypersurface.

We study the properties of such representations for singular (possibly reducible or non-reduced) hypersurfaces.
In particular, we obtain the decomposability criteria for determinantal representations of globally reducible
hypersurfaces.

Further, we classify the determinantal representations in terms of the corresponding kernel sheaves on $X$.
Finally, we extend the results to the case of symmetric/self-adjoint representations, with
implications to hyperbolic polynomials and generalized Lax conjecture.
\end{abstract}
\tableofcontents
\section{Introduction}\label{Sec Introduction}
\subsection{Setup}
Let $\k$ be an algebraically closed, normed, complete field of zero characteristic, e.g. the complex numbers, $\C$.
 Let $\k^n$ be the corresponding affine space, let $(\k^n,0)$ be the germ at the origin, i.e. a small neighbourhood.
 Let $\cOn$ denote the corresponding local ring of regular functions, i.e.
\\ -  rational functions that are regular at the origin, $\k[x_1,..,x_n]_{(\cm)}$, or
\\ - locally converging series, $\k\{x_1,..,x_n\}$, or
\\ - formal series $\k[[x_1,..,x_n]]$.

We denote the identity matrix by $\one$ and the zero matrix by $\zero$.
Let $\cM$ be a $d\times d$ matrix with the entries in either of:
\\ - (local case) $\cOn$
\\ - (global case) linear forms in $x_0,\dots,x_n$ (the later being the homogeneous coordinates of $\P^n$), i.e.
the global sections of the line bundle $\cO_{\P^n}(1)$.

We always assume $f:=\det(\cM)\not\equiv0$ and $d>1$. Such a matrix defines:
\\  - (local case) the germ of hypersurface near the origin, $(X,0):=\{\det(\cM)=0\}\sset(\k^n,0)$,
\\ -  (global case) the projective hypersurface $X:=\{\det(\cM)=0\}\sset\P^n$.

This hypersurface is called determinantal and the matrix $\cM$ is its \dr. The determinant, $f=\det\cM$, can be reducible
or non-reduced (i.e. not square-free). Let $f=\prod f^{p_\al}_\al$ be the (local/global) decomposition, i.e. $\{f_\al\}$
 are reduced, irreducible and mutually prime. Correspondingly the hypersurface is (locally/globally) decomposable:
$(X,0)=\cup (p_\al X_\al,0)\sset(\k^n,0)$ or $X=\cup p_\al X_\al\sset\P^n$.
Sometimes we consider the {\em reduced locus}: $X_{red}=\cup X_\al=\{\prod f_\al=0\}$.

The local/global \dr s are considered up to the {\em local/global equivalence}: $\cM\sim A\cM B$,
where $A,B\in GL(d,\cOn)$ or $A,B\in GL(d,\k)$. Both equivalences preserve the hypersurface pointwise.
\\
\\

Such "matrices of functions" appear constantly in various fields. Hence the interest in the \dr s of {\em arbitrary}
hypersurfaces (not only smooth or irreducible).
In this work we study the global \dr s of singular (possibly reducible, non-reduced) hypersurfaces/plane curves.
The symmetric and self-adjoint \dr s are treated separately at the end of the paper.
\subsection{A brief history}\label{Sec Intro Brief Overview}
A good summary of 19'th century's works on \dr s is in \cite{Wall78}. A modern introduction is in \cite{Dolgachev-book}.
\li
The question "For which pairs $(n,d)$ is the generic hypersurface of degree $d$ in $\P^n$ determinantal?"
has been studied classically. Already \cite{Dickson21} has shown that this happens only for $(2,d)$ and $(3,d\le3)$.
For a recent description see \cite{Beauville00}, \cite[\S4]{Dolgachev-book}  and
 \cite{Kosir2003}.
\li
Any (projective) curve in $\P^2$ admits a {\em symmetric} determinantal representation.
For smooth curves this was constructed (using ineffective theta characteristics) in \cite{Dixon1900}.
For singular curves this was proved in \cite[\S2 and \S7]{Barth77} and in \cite[prop. 2.28]{Catanese81}.
For some related works  see \cite{Room38}, \cite{Arbarello-Sernesi1979}.

For {\em smooth} (irreducible, reduced) plane curves the ordinary/symmetric/self-adjoint
 \dr s have been classified in \cite{Vinnikov89}\cite{Vinnikov93},
 see also \cite[prop.1.11 and cor.1.12]{Beauville00} and \cite[\S4]{Dolgachev-book}).
In \cite[theorem 3.2]{Ball-Vinni96} the classification of ordinary \dr s was extended to
the case of multiple nodal curves
 i.e. curves of the form $\{f^p=0\}\sset\P^2$ for $p\in\mN$ and $\{f=0\}$
irreducible, reduced, nodal curve.
\li A cubic surface in $\P^3$ is determinantal iff it contains at least two lines \cite[Proposition 4.3]{Brundu-Logar1998}.
In particular, the only cubic surfaces not admitting \dr s are those with a singularity of $E_6$ type, e.g.
 $\{x_0x^2_1+x_1x^2_2+x^3_3=0\}\sset\P^3$. For the classification of \dr s of smooth cubics  cf.
\cite{Buckley-Košir2007}, in \cite[\S9.3]{Dolgachev-book} the classification was extended to all cubic surfaces.
\li Determinantal quartic surfaces in $\P^3$ form a subvariety of codimension one in the family of all
the quartics (i.e. the complete linear family $|\cO_{\P^3}(4)|$). Such a surface may have on it any number
of lines up to 64, \cite{Room1950}.
In \cite{Giacobazzi1997} one studies \dr s of quartics in $\P^3$
possessing two lines $L_1,L_2$ of multiplicities $mult_1+mult_2=4$.
\li
In higher dimensions the determinantal hypersurfaces are necessarily singular and
 the singular locus is of dimension at least $(n-4)$. (For symmetric \dr s the dimension is at least $(n-3)$,
 in fact the subset of $X$ over which the corank of $\cM$ is at least two is of dimension at least $(n-3)$.)
The singularities occurring
at the points of $corank\cM\ge2$ are called {\em essential}, all the others: {\em accidental}.
The general linear symmetric determinantal hypersurface of
degree $d$ has only essential singularities, \cite[pg.495]{Salmon1865}.
For their properties and the classification of singularities of determinantal cubic/quartic surfaces,
i.e. $n=3$, see \cite{Piontkowski2006}.
 Nodal quartics in $\P^4$ were studied in \cite{Pettersen-1998}.
\li
The symmetric \dr s can be considered as  $n$-dimensional linear families of quadrics in $\P^{d-1}$.
 Hence various applications to Hilbert schemes of complete intersections, see \cite{Tjurin75-lectures}.
In \cite{Wall78} the \dr s of plane quartics (corresponding to nets of quadrics in $\P^3$) are studied in details.
\li
 The natural objects associated to a \dr\ are the kernel and cokernel of $\cM$. At each point of $X$ these are just
 vector spaces, as the point travels along the hypersurface these spaces glue into torsion-free sheaves supported
 on on $X$.
The \dr\  is determined (up to  the local/global equivalence) by its kernel/cokernel, e.g.
\cite[Thm 1.1]{Cook-Thomas79}. For the precise definition see \S\ref{Sec.Preliminaries.Kernel.Sheaves}.

\li As was proved in \cite{Helton-McCullough-Vinnikov2006}, any affine hypersurface in $\k^n$ admits a
symmetric \dr, i.e. any polynomial $f(x_1,..,x_n)$
can be presented as the determinant of a symmetric matrix of the type $A_0+\sum x_i A_i$.
\li
There are several reasons to consider non-reduced hypersurfaces, i.e. the cases when $\det\cM$ is not square-free.
For example, consider {\em matrix factorizations}, \cite{Eisenbud80}:
$AB=f\one$ with $\det(A)=$(a power of $f$). So matrix factorizations correspond to some \dr s of hypersurfaces
 with multiple components.
 And while the general hypersurface in $\P^n$ does not admit a \dr, unless $(d,n)=(3,3)$ or $n=2$,
its higher multiples, $\{f^p=0\}\subset\P^n$, do. For example, by \cite{BackHerzSand88},
\cite{Herzog-Ulrich-Backelin1991}, any homogeneous polynomial $f$ admits a matrix factorization in
linear matrices: $f\one=\cM_1\cdots\cM_d$, i.e. all the entries of $\{\cM_\al \}$ are linear.

\li The problem can be reformulated
as the study of $(n+1)$-tuples of matrices up to the two-sided equivalence, $(\cM_0,\dots,\cM_n)\sim A(\cM_0,\dots,\cM_n)B$.
Hence the applications in
linear algebra, operator theory and dynamical systems (see e.g. \cite{Ball-Vinni96}, \cite{Ball-Vinni2003},
 \cite{Tannenbaum81} or \cite{Livšic-Kravitsky-Markus-Vinnikov-book}).
In particular, these applications ask for the properties of determinantal representations of an {\em arbitrary}
hypersurface, i.e. with arbitrary singularities, possibly reducible and non-reduced.

\li  In applications one meets \dr s with specific properties, e.g. symmetric
or self-adjoint (in the real case). The self-adjoint \dr s are important in relation to the Lax conjecture as
they produce hyperbolic polynomials, see
\cite{Lax1958}, \cite{Guler1997},  \cite{Bauschke-Guler-Lewis-Sendov2001}, \cite{Lewis-Parrilo-Ramana2005},
\cite{Renegar2006}, \cite{Netzer-Thom2010}, \cite{Brändén2010}.

\li Finally, we mention the fast developing field of semi-definite-programming and matrix inequalities, i.e.
 presentability of the boundary of a convex set
in $\R^n$ by the determinant of a self-adjoint, positive definite matrix. For the introduction cf.
\cite{Helton-Vinnikov2007}, \cite{Livšic-Kravitsky-Markus-Vinnikov-book}.
\subsection{Results and Contents of the paper}
We tried to make the paper readable by non-specialists in commutative algebra/algebraic geometry. Thus in
\S\ref{Sec.Preliminaries.and.Notations} and further in the paper we recall some notions and results.
In particular in \S\ref{Sec.Preliminaries.Sheaves.on.Sing Curves} we recall sheaves
on singular (possibly reducible and non-reduced) hypersurfaces and Hirzebruch-Riemann-Roch theorem for
locally free sheaves.

In this paper we study the {\em global} \dr s. But the {\em local} version of the problem appears constantly,
due to the presence of singular points of curves/hypersurfaces and points where the kernel sheaves are not locally free.
The relevant results on the local version of the problem are obtained in \cite{Kerner-Vinnikov2010} and are
restated in \S\ref{Sec.Preliminaries.Local.det.reps}. Every global \dr\ can be localized and every local {\em algebraic}
\dr\ is comes from a global one. The localization process preserves the equivalence in a strong sense, etc.

In \S\ref{Sec.Preliminaries.Kernel.Sheaves} we introduce the sheaves of kernels (or kernel modules in the local case)
and prove some of their properties. The kernel sheaves can be also defined in a completely geometric way as follows.
Taking the kernel of a matrix provides a natural map $X\ni pt\to Ker(\cM|_{pt})\sset\k^d$. The image of $X$ in $\P^{d-1}$
under this map determines the kernel sheaf. This map is studied in \S\ref{Sec.Kernel.Sheaf.Geometric.Definition},
the equivalence of the two definitions is proven in proposition
\ref{Thm.Kernel.Geometric.Definition.Reduced.Case}.
Then we study particular types of \dr s/kernel sheaves: {\em \mg}\ \dr s (in \S\ref{Sec.Preliminaries.Maximaly.Generated.DR.s})
and {\em \XXS} (in \S\ref{Sec.Preliminaries.XXS.DR.s}). They possess especially nice properties and tend to
be {\em decomposable}.
\subsubsection{Decomposability.}
Suppose the determinant is reducible, $\det\cM=f_1f_2$, so the corresponding hypersurface is globally reducible:
$X=X_1\cup X_2$
Is $\cM$ globally decomposable? Namely, is it
globally equivalent to a block-diagonal matrix with blocks defining the components of the hypersurface:
$\cM\sim\cM_1\oplus\cM_2$. We study this question in \S\ref{Sec Global Decomposability}.
The global decomposability obviously implies the local one.
A probably unexpected feature is the converse implication:
\\
{\bf Theorem} ~ \ref{Thm Decomposability Global from Decomposability Local}
{\em Let $X=X_1\cup X_2\sset\P^n$ be a global decomposition of the hypersurface. Here $X_1,X_2$ can
be further reducible, non-reduced, but without common components, i.e their defining polynomials are
relatively prime.
 $\cM$ is globally decomposable, i.e. $\cM\stackrel{globally}{\sim}\cM_1\oplus\cM_2$,
 iff it is locally decomposable at each point $pt\in X_1\cap X_2$, i.e.
$\cM\stackrel{locally}{\sim}\one\oplus\cM_1|_{(\P^n,pt)}\oplus\cM_2|_{(\P^n,pt)}$.
Here $\cM_\al|_{(\P^n,pt)}$ is the local \dr\  of $(X_\al,pt)$.
}

The proof of this property is heavily based on Noether's $AF+BG$ theorem \cite[pg. 139]{ACGH-book},
in fact one might consider the statement as Noether-type theorem for matrices.

Similarly, suppose at each point of the intersection $X_1\cap X_2$ the \dr\ is
locally equivalent to an upper-block-triangular, $\cM\sim\bpm\cM_1& *\\\zero&\cM_2\epm$,
where $\det(\cM_\al)$ defines $X_\al$. Then the
global \dr\  is globally equivalent to an upper-block-triangular,
proposition \ref{Thm Decomposability Global From Sheaf Axioms}.

Both statements are non-trivial from linear algebra point of view, but almost tautological
when considered as statements on
kernel sheaves. In the first case the sheaf is the direct sum, $E\approx E_1\oplus E_2$,
in the second case it is an extension:  $0\to E_1\to E\to E_2\to0$.

These results completely reduce the global decomposability problem for \dr s to the local problem.
In \S\ref{Sec.Preliminaries.Maximaly.Generated.DR.s} and \S\ref{Sec.Preliminaries.XXS.DR.s} we state
various necessary and sufficient criteria for local decomposability of \dr s, they are formulated and proved
in \cite{Kerner-Vinnikov2010}.
\subsubsection{Properties and classification of kernel sheaves.}
In \S\ref{Sec.Kernels.Properties.Classification} we study the kernel sheaves on the hypersurfaces in $\P^n$.
First we summarize their properties Theorem \ref{Thm.Kernel.Sheaves.Properties}.
It is possible to classify those sheaves arising as kernels of \dr s.
For the smooth case this was done in \cite{Vinnikov89} see also \cite[\S4]{Dolgachev-book}.
The classification for an arbitrary hypersurface is done in \cite[Theorem A]{Beauville00}.
We give a direct proof of this result.
\\$\Big(${\bf Theorem} \ref{Thm.Kernel.Sheaves.Classification} and {\bf Theorem}
 \ref{Thm.Kernel.Sheaves.Properties}$\Big)$
{\em Consider a hypersurface $X=\cup p_\al X_\al\sset\P^n$. The torsion-free sheaf $E_X$ of multi-rank $(p_1,..,p_k)$
  is the kernel of a \dr\ of $X$, in the sense of equation
 \eqref{Eq.Kernel.Sheaf.Algebraic.Definition}, iff $h^0(E_X(-1))=0$, $h^i(E_X(j))=0$ for $0<i<n-1$, $j\in\Z$
 and $h^{n-1}(E_X(1-n))=0$.}
\\(Note that over a non-reduced hypersurface the torsion-free sheaf can be nowhere locally free. For the definition
of multi-rank see \S\ref{Sec.Preliminaries.Rank.of.Sheaf})

We prove this theorem by an explicit construction, generalizing \cite{Dixon1900} and \cite{Vinnikov89}, where it
is done for smooth plane curves.
One advantage of this proof is that it is easily adjustable to symmetric/self-adjoint cases (see below).

An immediate corollary, in the case when $E_X$ is locally free, is the information about the Chern class of the kernel,
theorem \ref{Thm.Kernel.Sheaves.Properties}.
\subsubsection{Relation to matrix factorizations and descent to the reduced locus.}
Suppose $AB=\one\prod f_\al$, where $A$ has homogeneous linear entries, the factors $\{f_\al\}$ are irreducible
and $A$ is non-invertible
at any point of the hypersurface $\{f=0\}$. Then $\det(A)=\prod f_\al^{p_\al}$, for some multiplicities $\{p_\al\}$.
So, $A$ is a \dr\ of the non-reduced hypersurface. A natural question: which \dr s arise from matrix factorizations
of reduced hypersurfaces?
An immediate consequence of our results:
\bcor Let $\cM$ be a \dr\ of $\prod f_\al^{p_\al}$. There exists a matrix $\cN$ satisfying
$\cM\cN=\prod f_\al\one$ iff $\cM$ is \mg\ at generic smooth points of the reduced locus $\{\prod f_\al=0\}$.
\ecor
(For the definition of \mg\ see \S\ref{Sec.Preliminaries.Maximaly.Generated.DR.s}.)
In such a case the kernel  $E$ of $\cM$, a sheaf over the non-reduced hypersurface $X$, has natural descent to
the reduced locus $X_{red}$: $E\rightsquigarrow E^{red}_X$. In \S\ref{Sec.Kernel.Sheaf.On.Reduced.Locus} we
classify the sheaves obtained in this way:
\\{\bf Proposition \ref{Thm.Kernel.Sheaves.Reduced.Locus.Classification}} {\em A torsion free sheaf $E^{red}_X$ of
multirank $(p_1,..,p_k)$ on the reduced locus $X_{red}$ arises by descent from $X$ iff
$h^0(E^{red}_X(-1))=0$, $h^i(E^{red}_X(j))=0$ for $0<i<n-1$, $j\in\Z$ and $h^{n-1}(E^{red}_X(1-n))=0$.}
\subsubsection{Ascent to the modification, for curves}
Let $C=\cup p_\al C_\al$ be the global decomposition of a plane curve. One often considers normalization:
$\tC:=\coprod (p_\al\tC_\al)\norm\cup p_\al C_\al$, here each $\tC_\al\to C_\al$ is the normalization of an irreducible curve.
  Correspondingly the kernel sheaf is pulled back: $\nu^*(E)/Torsion$.
For the normalization $\tC\norm C$ the pullback $\nu^*(E)/Torsion$
is locally quasi-free (or just locally free in the reduced case).

Sometimes the pullback is locally quasi-free already for some intermediate modification: $\tC\to C'\norm C$.
It is important to classify those sheaves on $C'$ whose pushforward to $C$ produces kernel sheaves.
\\{\bf Corollary \ref{Thm.Kernel.Descent.From.Modification.Classification}} {\em Given a modification
$C'\to\cup p_\al C_\al$, the torsion free sheaf $E_{C'}$  descends to the kernel of
a \dr\ of $C$ iff $h^i(E^{red}_{C'}(-1))=0$ for $i\ge0$.}

A more complicated question is: which sheaves on $C'$ are pullbacks (modulo torsion) of kernel sheaves on $C$?
(Note that in general $E\subsetneq\nu_*\nu^*(E)/Torsion$.)
We give a criterion in proposition \ref{Thm.Kernel.Ascent.To.Modification}.
\\\

Once the general properties of kernel sheaves are established one can study the \dr s for particular hypersurfaces.
In \S\ref{Sec.Kernel.Sheaves.Examples} we give some simplest examples of kernel sheaves on curves/surfaces.

\subsubsection{Symmetric and self-adjoint \dr s}
In \S\ref{Sec.Symmetric.DetReps} and \S\ref{Sec.Self.Adjoint.dr s} we work with $\k=\R\sset\C$.
If $\cM$ is symmetric or self-adjoint then it is natural to consider symmetric or
self-adjoint equivalence ($\cM\stackrel{s}{\sim}A\cM A^T$ or $\cM\stackrel{\tau}{\sim}A\cM A^\tau$).
Many of the previous results are extended to this setup.

 Being symmetric or self-adjoint can be formulated  in terms of the kernel sheaves
 (properties \ref{Thm.Kernel.is.Symmetric.iff.Det.Rep.Is.Symmetric} and
  \ref{Thm.Self-adjoint.Equivalence.Iff.Ordinary.Equiv}).
Two symmetric representations are equivalent (in the ordinary sense) iff they are symmetrically equivalent (proposition
\ref{Thm.Symmetric.Equivalence.iff.Ordinary.Equivalence}). For self-adjoint representations this is true up
to a diagonal matrix, the precise statement is proposition \ref{Thm.Self-adjoint.Equivalence.Iff.Ordinary.Equiv}.

The symmetric \dr s of singular hypersurfaces are studied in \S\ref{Sec.Symmetric.DetReps}.
In particular, we characterize the kernel sheaves of symmetric \dr s of hypersurfaces.
In \S\ref{Sec.Self.Adjoint.dr s} we characterize the self-adjoint \dr s of hypersurfaces.

\subsubsection{Applications to hyperbolic polynomials}
Recall that if $\cM$ is self-adjoint then $\det\cM$ is a hyperbolic polynomial,
\S\ref{Sec.Self Adjoint dr of hyperbolic polynomials}.
Then the real locus of $X$ can have
at most one singular point with a non-smooth locally irreducible component,
theorem \ref{Thm.Hyperbol.Curve.Has.Smooth.Branches.Det.Rep.is.Maximal}. In the later case the
region of hyperbolicity
degenerates to this singular point. Thus, if the hypersurface is defined by a self-adjoint positive-definite \dr,
then all the locally irreducible components of its reduced locus are smooth.
In theorem \ref{Thm.Hyperbol.Curve.Has.Smooth.Branches.Det.Rep.is.Maximal} we prove that any self-adjoint
positive-definite \dr\  of a real hypersurface is $\tX/X$ saturated (at real points),
i.e. its kernel arises as the push-forward of a locally free sheaf from the normalization $\tX\norm X$.

\section{Preliminaries and notations}\label{Sec.Preliminaries.and.Notations}
For local considerations we always assume the (singular) point to be
at the origin and mostly use the ring of locally convergent power series
$\k\{x_1,..,x_n\}=\cOn$. Let $\cm=\bl x_1,\dots,x_n\br\sset\cOn$ be the maximal ideal.

The {\em tangent cone} $T_{(X,0)}\sset\k^n$ is formed as the limit of all the
tangent hyperplanes at smooth points. For the hypersurface $(X,0)=\{f=0\}$, with the Taylor expansion
 $f=f_p+f_{p+1}+\cdots$, the tangent cone is $\{f_p=0\}\sset(\k^n,0)$.
For curves the tangent cone is the collection of tangent lines,
each with the corresponding multiplicity.

The tangent cone is in general reducible. Associated to it is the {\em tangential} decomposition:
$(X,0)=\cupl_{\al\in T_{(X,0)}}(X_\al,0)$. Here $\al$ runs over all the (set-theoretical) components of the tangent cone,
each $(X_\al,0)$ can be further reducible, non-reduced.
\bex
Consider the curve singularity $(X,0)=\{(y^2-x^4)(x^2-y^4)=0\}\sset(\k^2,0)$. Here the tangent cone
is $T_{(X,0)}=\{y^2x^2=0\}\sset\k^2$. Accordingly, the tangential decomposition is:
$\{y^2=x^4\}\cup\{x^2=y^4\}\sset(\k^2,0)$.
\eex

The basic invariant of the hypersurface singularity $\{f_p+f_{p+1}+\cdots=0\}$ is the {\em multiplicity} $mult(X,0)=p$.
 For the tangential components denote $p_\al=mult(X_\al,0)$.

\subsection{Sheaves on singular hypersurfaces}\label{Sec.Preliminaries.Sheaves.on.Sing Curves}
The theory of coherent sheaves on multiple smooth curves,
i.e. $pC_{red}$, for $C_{red}$ irreducible and smooth, is developed in \cite{Drezet2009}.

A coherent sheaf on a pure dimensional scheme $X$ is called torsion-free if it is has no subsheaf
whose support is of strictly smaller dimension than $dim(X)$.

\subsubsection{Multi-rank of pure sheaves on reducible, non-reduced hypersurfaces}\label{Sec.Preliminaries.Rank.of.Sheaf}
Let $F_X$ be a torsion-free sheaf, its singular locus $Sing(F_X)\subset X$ is the set of (closed) points
where $F_X$ is not locally free. If $X$ is reduced then $F$ is generically locally free
and to the decomposition $X=\cup X_\al$ is associated the multi-rank $(r_1,..,r_k)$: $r_i=rank(F|_{X_i})$.

In the non-reduced case a torsion free sheaf can be {\em nowhere locally free}. To define its multi-rank we need a
preliminary construction. Consider a multiple hypersurface, $X=p X_{red}=\{f^p=0\}$, where $X_{red}$
is irreducible. We define the rank of $F_X$. Let $I_{X_{red}}\sset\cO_X$ be the ideal of the reduced locus.
As $X$ is a hypersurface, this ideal is principal, $I_{X_{red}}=\bl f\br$. Consider the multiplication by $f$ on $F$.
Its successive kernels define a useful filtration on $F$:
\beq
0\subseteq Ker(f)\subseteq Ker(f^2)\subseteq\cdots\subseteq Ker(f^{p-2})\subseteq Ker(f^{p-1})\subseteq Ker(f^p)=F
\eeq
Associated to this filtration is the graded sheaf:
\beq
Gr_f(F):=\oplus^{p-1}_{j=0}\quotient{Ker(f^{p-j})}{Ker(f^{p-j-1})}=\oplus^{p-1}_{j=0} Gr_j
\eeq
By definition $fGr_j=0$, hence $Gr_j$ is naturally a module over $\cO_{X_{red}}$.
\bex
$\bullet$ Though $F$ is torsion-free, its reduction, $F\underset{\cO_X}{\otimes}\cO_{X_{red}}=F/IF$, in general has torsion. For example,
let $X=\{x_1^2=0\}\sset\P^n$. Let $F$ be an $\cO_X$ module, generated by $s_1=\bpm x_1\\0\epm$ and $s_2=\bpm x_2\\x_1\epm$.
So $F=\cO_X\bl s_1,s_2\br/(x_1s_2-x_2s_1,x_1s_1)$. Then $F/IF=\cO_{X_{red}}\bl s_1,s_2\br/(x_2s_1)$, i.e. the
element $s_1$ is annihilated by a non-zero divisor $x_2$.
\li For simplicity consider the local case of curves. Let $\cO_C=\k[x,\ep]/\ep^p$
and $F=\oplusl^{p-1}_0 \ep^j(\k[x])^{\oplus l_j}$, for $l_0\le l_1\le .. l_{p-1}$. Then the filtration is:
\beq
0\sset\ep^{p-1}(\k[x])^{\oplus l_{p-1}}\sset\cdots\sset\oplusl^{p-1}_1 \ep^j(\k[x])^{\oplus l_j}\sset
\oplusl^{p-1}_0 \ep^j(\k[x])^{\oplus l_j}
\eeq
\eex
A sheaf is called {\em locally quasi-free} (at a point) if its graded version is locally free (at this point).
Every torsion free sheaf on a hypersurface is generically locally quasi-free.
\bed
1. Let $X=pX_{red}$, where $X_{red}$ is irreducible and reduced.
The rank of $F$ on $X$ is the rank of $Gr_f(F)$ as a module over $X_{red}$.
\\2. For $X=\cup p_\al X_\al$ the multi-rank of $F$ is the collection of ranks $\{rank F|_{p_\al X_\al}/Torsion\}$.
\eed
\bprop\label{Thm.Rank.of.Sheaf.Equals.Intersection.Multiplicity}
Let $L\sset\P^n$ be the generic line, let $pt\in X\cap L$. Then $length(F|_{pt})=rank(F)$.
\eprop

\subsubsection{Hirzebruch-Riemann-Roch theorem}
In this paper we use Hirzebruch-Riemann-Roch theorem for locally free sheaves on (singular, reducible, possibly
non-reduced) hypersurfaces, \cite[pg.354]{Fulton-book}:
\beq
\chi(F_X)=\Big(ch(F_X)Td(T_X)\Big)_{top.dimensional}
\eeq
Here
\ls the Euler characteristic of the sheaf is $\chi(F_X)=\sum^{n-1}_{i=0}(-1)^ih^i(F_X)$.
\ls the Chern character of the sheaf is $ch(F_X)=\prod\exp(\al_i)$ where $\{\al_i\}_i$ are the Chern roots of $F_X$,
i.e. $c(F_X)=\prod(1+\al_i)$.
\ls the Todd class of the hypersurface $Td(T_X)$. For a smooth variety it equals
$\prod_i\frac{\al_i}{1-\exp(-\al_i)}$, where $\{\al_i\}_i$ are the Chern roots of the tangent bundle, i.e.
$c(T_X)=\prod(1+\al_i)$.
Suppose the scheme $X$ is singular but embeddable as a locally complete intersection into a smooth
variety, $X\sset Y$.  For example, this is the case for hypersurfaces in $\P^n$.
Then $Td(T_X)$ is the Todd class of the virtual tangent bundle, $T_X:=T_Y|_X-N_{X/Y}$.
\ls both the Chern and the Todd classes are graded, from their product one extracts the top dimensional part.
\bex\label{Ex.Chern.Class.Hypersurface}
Let $X_d\sset\P^n$ be an arbitrary hypersurface of degree $d$. Let $\cL$ be a line bundle on $X$. The
total Chern class is $c(\cL)=1+c_1(\cL)$. Then
$ch(\cL)=\sum \frac{c^i_1(\cL)}{i!}$. The virtual tangent bundle of $X$ is defined by
\beq
0\to T_X\to T_{\P^n}|_X\to N_{X/\P^n}\to0,\qquad c(T_{\P^n})=(1+L)^{n+1},\qquad c(N_{X/\P^n})=1+dL
\eeq
Here $L$ is the class of a hyperplane in $\P^n$. Hence
\beq
c(T_X)=\frac{(1+L)^{n+1}}{1+dL}=1+(n+1-d)L+\Big(\bin{n+1}{2}-d(n+1)+d^2\Big)L^2+..
\eeq
The Hirzebruch-Riemann-Roch theorem in this case reads:
\beq
\chi(\cL)=\Bigg[\big(\sum \frac{c^i_1(\cL)}{i!}\big)\big(1+\frac{c_1(T_X)}{2}+\frac{c^2_1(T_X)+c_2(T_X)}{12}+
\frac{c_1(T_X)c_2(T_X)}{24}+..\big)\Bigg]_{top.dim.}
\eeq
For example, in the case of plane curves, \cite[\S IV.I  exercise 1.9]{Hartshorne-book}:
\beq\label{Eq Riemann Roch for Coh Sheaves}
h^0(F)-h^1(F)=deg(F)+(1-p_a)rank(F)
\eeq
Here $p_a=\binom{d-1}{2}$ is the arithmetic genus of the plane curve, it does not depend on the singularities.
\eex

The same formula holds sometimes for sheaves that are torsion free, but not locally free.
For example, for torsion-free sheaves on an integral curve this was proved in \cite[thm 1.3]{Hartshorne86}.
See also \cite{Fulton-02}.
 The theorem was also proved for sheaves on multiple smooth curves in \cite{Drezet2009}.
\subsubsection{The dualizing sheaf and Serre duality}
For torsion free sheaves on varieties with (at most) Gorenstein singularities, e.g. on any hypersurface in $\P^n$,
 the dualizing sheaf $w_C$ is invertible. By the adjunction formula for a hypersurface in $\P^n$
 of degree $d$, with arbitrary singularities: $w_X=\cO_X(d-n-1)$.
Then the usual Serre duality holds: $H^i(F_X)=H^{\dim(X)-i}(F^*_X\otimes w_X)^*=H^{\dim(X)-i}\big(F^*_X(d-n-1)\big)^*$.

\subsection{The matrix and its adjoint}\label{Sec.Preliminaries.Matrix.and.its.Adjoint}
We work with (square) matrices, their sub-blocks and particular entries. Sometimes to avoid confusion we
 emphasize the dimensionality, e.g. $\cM_{d\times d}$. Then $\cM_{i\times i}$ denotes an $i\times i$
block in $\cM_{d\times d}$ and $\det(\cM_{i\times i})$ the corresponding minor. On the other hand by $\cM_{ij}$
we mean a particular entry.
\\
\\
Let $\cM$ be a determinantal representation of $X\sset\P^n$ or  $X\sset\k^n$. Let $\cMv$ be the adjoint matrix of $\cM$,
so $\cM\cMv=\det(M)\one_{d\times d}$.
Then $\cM$ is non-degenerate outside the hypersurface $X$ and its corank over the hypersurface satisfies:
\beq\label{Eq Bound on Corank of Matrix}
1\le corank(\cM|_{pt\in X})\le mult(X,pt)
\eeq
(as is checked e.g. by taking derivatives of the determinant).
The adjoint matrix $\cMv$ is not zero at smooth points of $X$.
As $\cMv|_X\times\cM |_X=\zero$ the rank of $\cMv$ at any smooth point of $X$ is 1 (for the reduced hypersurface).
Note that $(\cM ^\vee)^\vee=f^{d-2}\cM $ and $\det\cMv=f^{d-1}$.

\subsubsection{The case $det(\cM)\equiv0$} A natural question in this case if whether $\cM$ is equivalent to a matrix with a zero
row/column. In general this does not hold, e.g. for $\cM=\bpm 0&0&x\\0&0&y\\x&z&0\epm$. Indeed, if $\cM$ was equivalent
to a matrix with zero row/column then $\cMv$ would be equivalent to a matrix with at most 3 non-zero entries.
But the ideal of $\cMv$ is $\bl x^2,xy,xz,yz\br$, i.e. is generated by 4 elements. And this ideal is invariant
under equivalence.
\subsection{Local \dr s}\label{Sec.Preliminaries.Local.det.reps}
Here we review some aspects of local \dr s and quote the necessary results, all the proofs are
in \cite{Kerner-Vinnikov2010}.
Essentially this is the part of commutative algebra, the theory of \CM modules,
see \cite{Yoshino-book} and \cite{Leuschke-Wiegand-book}.

\subsubsection{The global-to-local reduction}\label{Sec.Preliminaries.Local.vs.Global}
This is the way to pass from global to local \dr s.
Replace the homogeneous coordinates of $\P^n$ by the local coordinates: $(x_0,..,x_n)\to(x_1,..,x_n)$ with $x_0=1$.
\bpro\label{Thm Localiz Chip off Unity}
Suppose the multiplicity of $(X,0)$ is $m\ge1$ and $\cM_{d\times d}$ is a corresponding (local or global) determinantal
representation.
\\1. Locally $\cM_{d\times d}$ is equivalent to $\bpm \one_{(d-p)\times(d-p)}&\zero\\\zero& \cM_{p\times p}\epm$
with $\cM_{p\times p}|_{(0,0)}=\zero$ and $1\le p\le m$.
\\2. The stable local equivalence (i.e. $\one\oplus\cM_1\sim\one\oplus\cM_2$) implies ordinary local
equivalence ($\cM_1\sim\cM_2$). So, the global-to-local reduction is unique up to the local equivalence.
\epro
From the algebraic point of view the first statement is the reduction to the minimal free resolution of the
kernel module \cite[\S20]{Eisenbud-book}.
 The first claim is proved in symmetric case e.g. in \cite[lemma 1.7]{Piontkowski2006}.
 Both bounds are sharp, regardless of the singularity of hypersurface. For the second statement
 see \cite{Kerner-Vinnikov2010}.
\bed\label{Def Localication}
In the notations as above, $\cM_{p\times p}$ is the reduction of $\cM_{d\times d}$ or the local representation.
\eed
Any matrix whose entries are {\em rational} functions, regular at the origin, is the
 reduction of some global determinantal representation:
\bel\label{Thm Localiz Every Local Arises From Global}
1.For any $\cM_{local}\in Mat(p\times p,\k[x_1,..,x_n]_{(\cm)})$, there exists a matrix of homogeneous linear forms,
$\cM_{global}\in Mat\Big(d\times d,H^0(\cO_{\P^N}(1))\Big)$, whose reduction is $\cM_{local}$.
\\2. In particular, if $\cM_1$, $\cM_2$ are locally equivalent and $\cM_1$ is the reduction of some $\cM_{global}$ then $\cM_2$
is also the reduction of $\cM_{global}$.
\eel
Note that in the lemma $\cM_{global}$ or $\det(\cM_{global})$ are not unique in any sense, even
the dimension of $\cM_{global}$ is not fixed.
\bpr We can assume that $\cM$ is a matrix with polynomial entries. Indeed, all the denominators of entries of $\cM$ do
not vanish at the origin, hence one can multiply $\cM$ by them.

Let $x^{a_1}_1..x_n^{a_n}$ be a monomial in $\cM_{local}$ with the highest total degree $\sum a_i$.
By permutation assume it belongs
to the entry $\cM_{11}$. Consider the augmented matrix:
\beq
\bpm1&0\\0&x^{a_1}_1..x_n^{a_n}+..&\cM_{12}&..\\0&\cM_{21}&\cM_{22}&..\\0&..&..\epm
\eeq
It is locally equivalent to
\beq
\bpm 1&x_1&0\\-x^{a_1-1}_1..x_n^{a_n}&0+..&\cM_{12}&..\\0&\cM_{21}&\cM_{22}&..\\0&..&..\epm
\eeq
For the new matrix the number of monomials with highest total degree is less by one. Continue in the
same way till all the monomials of the highest total degree ($\sum a_i$) are removed. Now the highest order degree is less than $\sum a_i$.
Continue by induction till one gets a matrix with entries of degree at most 1.
\epr
The last lemma is formulated as a purely linear-algebraic statement. A reformulation in terms of sheaves (using proposition
 \ref{Thm.Notions.of.equality of kernels}):
\\{\em For any kernel sheaf the stalk (at any point) is a kernel module. The isomorphism class of the stalk is well defined.
 Every kernel module is the stalk of some kernel sheaf.}

\bex\label{Ex.Det.Reps.of.A_k}
Let $\cM$ be a \dr\ of the plane curve $\{y^2=x^{k+1}\}\sset\k^2$, this is the $A_k$ singularity. Suppose $\cM$ is local,
i.e. $\cM|_{(0,0)}=\zero$. As the multiplicity of this curve singularity is 2, the dimensionality of $\cM$ is either 1 or 2.
The first case is trivial, in the second case one can show that $\cM$ is (locally) equivalent to $\bpm y&x^l\\x^{k+1-l}&y\epm$.
\eex
Finally we prove that for global \dr s the local equivalence is not weaker than the global one.
\bel
Suppose two global \dr s are locally equivalent, i.e. $\cM_1=A\cM_2 B$ for $A,B\in GL(d,\k[[x_0,\dots,x_n]])$.
Then $\cM_1$, $\cM_2$ are globally equivalent too.
\eel
\bpr
Expand $A=jet_0(A)+A_{\ge1}$, where $A_{\ge1}|_0=\zero$, similarly $B=jet_0(B)+B_{\ge1}$. Note that
$jet_0(A),jet_0(B)\in GL(d,\k)$. Therefore
\beq
\Big(jet_0(A)\Big)^{-1}\cM_1\Big(jet_0(B)\Big)^{-1}=(\one+A'_{\ge1})\cM_2(\one+B'_{\ge1})
\eeq
Hence, by comparing the degrees (in $x_i$) we get: $\Big(jet_0(A)\Big)^{-1}\cM_1\Big(jet_0(B)\Big)^{-1}=\cM_2$.
\epr
\subsection{Kernels and cokernels of \dr s}\label{Sec.Preliminaries.Kernel.Sheaves}
Let $\cM_{d\times d}$ be a \dr\  of the hypersurface $X\sset\P^n$. At each point of $X$ the matrix has some (co-)kernel.
These vector spaces glue to sheaves on $X$, or to vector bundles in nice situations. The sheaf structure can
be defined in two equivalent ways.

\subsubsection{Algebraic definition of the kernel}\label{Sec.Kernel.Sheaf.Algebraic.Definition}
The cokernel sheaf is defined by the sequence
\beq\label{Eq Exact Sequence for Sheaves on P2}
0\to\cO^{\oplus d}_{\P^n}(-1)\stackrel{\cM}{\to}\cO^{\oplus d}_{\P^n}\to Coker\to0
\eeq
As $\cM$ is invertible at the points of $\P^n\smin X$ the cokernel is supported on the hypersurface.
Restrict the sequence to the hypersurface (and twist),
then the kernel appears.
\beq\label{Eq.Kernel.Sheaf.Algebraic.Definition}
0\to E_X \to\cO_X^{\oplus d}(d-1)\stackrel{\cM}{\to}\cO_X^{\oplus d}(d)\to Coker(\cM)_X\to0
\eeq
Sometimes we consider also  the "left" kernel, $E_X^l$, the kernel of $\cM^T$, (called the Auslander transpose):
\beq
0\to E_X^l \to\cO_X^{\oplus d}(d-1)\stackrel{\cM^T}{\to}\cO_X^{\oplus d}(d)\to Coker(\cM^T)_X\to0
\eeq
From now on all the sheaves are considered on curves/hypersurfaces.
\bex\label{Ex.Det.Reps.Smooth.Quadrics.P3}
Consider a smooth quadric surface $X=\{x_0x_1=x_2x_3\}\sset\P^3$. By direct check, it has two (non-equivalent)
\dr s: $\bpm x_0&x_2\\x_3&x_1\epm$ and $\bpm x_0&x_3\\x_2&x_1\epm$. Consider the first case, the kernel $E_X$ is the line
bundle spanned by two sections: $\bpm -x_3\\x_0\epm$ and $\bpm -x_1\\x_2\epm$. To identify this line bundle
recall that $X\approx\P^1_{left}\times\P^1_{right}$ and the isomorpism can be written explicitly:
$(x_0,x_1,x_2,x_3)\to\Big((x_0,x_2),(x_0,x_3)\Big)$. Note that both maps are well defined, using
$(x_0,x_2)=(x_3,x_1)$ and $(x_0,x_3)=(x_2,x_1)$. The sections of $E_X$ vanish at $x_0=0=x_3$ and $x_2=0=x_1$.
Note that both cases define the divisors $pt\times\P^1_{right}\sset X$. Thus:
$E_X\approx \cO_{\P^1_{left}}(1)\boxtimes\cO_{\P^1_{right}}$.
\eex

To work with singular points we consider the stalks of the kernel sheaves, i.e. kernel modules over the local ring $\cO_{(X,0)}$.
For them one has the corresponding exact sequence of modules. Many properties of kernels hold
both in local and in global situation, usually we formulate and prove them together.

Both in local and global cases the kernel module/sheaf is spanned by the columns of the adjoint matrix $\cMv$, see
theorem \ref{Thm.Kernel.Sheaves.Properties}.

\subsubsection{Geometric definition of the kernel}\label{Sec.Kernel.Sheaf.Geometric.Definition}
The kernel sheaves $E_X,E^l_X$ can be defined also in a more geometric way \cite[\S3]{Vinnikov89}.
Suppose $X$ is reduced, so for generic point $pt\in X$ the kernel $Ker(\cM|_{pt})$ is a one-dimensional vector subspace
of $\k^d$.
 Consider the rational map $\phi:X\dashrightarrow\P^{d-1}$ defined on the smooth points of $X$ by:
\beq
X\ni pt\stackrel{\phi}{\to}\{Ker(\cM|_{pt})\sset\k^d\}\to\P^{d-1}
\eeq
It extends to a morphism of algebraic varieties iff $corank(Ker\cM|_X)\equiv1$, i.e. $E|_X$ is a
\parbox{11cm}{locally free sheaf.
In general, consider a birational morphism resolving the singularities of the map $\phi$, see the diagram.
Hence $\nu^*(E)/Torsion$ is a locally free sheaf and the pull-back $\nu^*(\phi)$ extends to a morphism.
Recall that $\cO_{\P^{d-1}}(-1)$
is the tautological bundle, its fibre over $(x_0,\dots,x_{d-1})$ is
spanned by the vector $(x_0,\dots,x_{d-1})$.}
\hspace{0.5cm}$\bM \tX\hspace{0.8cm} \nu^*(E)/Torsion \hspace{0.2cm} \ \\\nu\downarrow\searrow \ \hspace{3.2cm} \ \\
X\stackrel{\phi}{\dashrightarrow}\P^{d-1}\ \ \ \cO_{\P^{d-1}}(-1)\eM$
\bprop\label{Thm.Kernel.Geometric.Definition.Reduced.Case}
 Let $X\sset\P^n$ be a reduced hypersurface.
\\1. If $E_X$ is locally free and $\phi:X\to\P^{d-1}$ is the corresponding morphism,
then $\phi^*(\cO_{\P^{d-1}}(-1))=E_X(1-d)\sset\cO^{\oplus d}_X$.
\\2. If $\nu^*(E_X)/Torsion$ is locally free and
$\tX\stackrel{\phi\circ\nu}{\to}\P^{d-1}$ is the corresponding morphism then $(\phi\circ\nu)^*(\cO_{\P^{d-1}}(-1))=\nu^*(E_X)(1-d)/Torsion$.
\\3. The kernel sheaf  $E_X$ is determined uniquely by $\phi^*(\cO_{\P^{d-1}}(-1))$.
\eprop
\bpr Use the definition, i.e. \eqref{Eq.Kernel.Sheaf.Algebraic.Definition}
1. In the locally free case $E_X(1-d)$ and $\phi^*(\cO_{\P^{d-1}}(-1))$ are two line subbundles
of $\cO^{\oplus d}_X$, whose fibres coincide at each point. So the bundles coincide tautologically.
\\
\\
2. Similarly, $\nu^*(E_X(1-d))/Torsion$ and $(\phi\circ\nu)^*(\cO_{\P^{d-1}}(-1))$ are line subbundles on
 $\cO^{\oplus d}_{\tX}$, with coinciding fibres.
\\\\
3. Suppose there are two kernel sheaves corresponding to $\phi^*(\cO_{\P^{d-1}}(-1))$, their restrictions
onto the smooth locus of $X$ coincide, as sub-sheaves of $\cO^{\oplus d}(d-1)|_{X\smin Sing(X)}$.
Then, by proposition \ref{Thm.Notions.of.equality of kernels}, part 1, we get
 two \dr s, $\cM_1$ and $\cM_2$, satisfying locally $\cM_1=A\cM_2$. Here the entries of $A$
are in $\cO_{(\k^n,0)\smin Sing(X)}$ and $A$ is locally invertible at each point of $(\k^n,0)\smin Sing(X)$.
So, each entry of $A$ is regular in codimension one, i.e. its possible locus of irregularity is of codimension
at least two. But then this entry, being a rational function, is regular on $(\k^n,0)$. Similarly, $\det(A)\neq0$
except possibly for a subset of codimension two. Thus $\det(A)\neq0$ on the whole $(\k^n,0)$.
Therefore the stalks of two kernel sheaves coincides everywhere on $X$.
\epr
For some \dr s of {\em non-reduced} hypersurfaces the kernel can be defined geometrically too, see
 \S\ref{Sec.Kernel.Sheaf.On.Reduced.Locus}.
\bex\label{Ex.Kernel.Maps.and.Images}
1. Let $X\sset\P^2$ be a line arrangement, consider the simplest \dr: the diagonal matrix $\cM=(l_1,\dots,l_d)$.
Here $\{l_i\}_i$ are linear forms defining the lines. The map $X\stackrel{\phi}{\dashrightarrow}\P^{d-1}$ is defined on the smooth locus of $X$,
it sends each line to a point in $\P^{d-1}$.
\\2. In general, for a determinantal hypersurface $X=\cup_\al X_\al\sset\P^n$, let
 $X^0\sset X$ be an open dense subset on which $\phi$ is defined. Then $\cM$ is decomposable,
 $\cM\sim\oplus\cM_\al$, iff $\overline{\phi(X^0)}=\coprod Span(X_\al)$, where the spans
 ($Span(X_\al)\sset\P^{d-1}$ is the minimal linear subspace that contains $\phi(X_\al)$)
 are mutually generic, i.e. $\forall \al:$ $Span(X_\al)\cap Span(\coprod'_{\be\neq\al}X_\be)=\empty$.
\eex

\subsubsection{Kernels vs \dr s}
Finally we formulate the relation between the embedded kernels and the \dr s.
 As always, in the local case, we assume $\cM|_0=\zero$.
\bprop\label{Thm.Notions.of.equality of kernels}
0. The kernel module of $\cM_{(X,0)}$ is an $\cO_{(X,0)}$ module minimally generated by the columns of the adjoint
matrix $\cMv_{(X,0)}$. Similarly, for the kernel sheaf of $\cM_X$, the columns of the adjoint matrix give the
natural basis of the space $H^0(E_X)$. In particular $h^0(E_X)=d$.
\\1. Let $\cM_1,\cM_2\in Mat(d\times d,\cOn)$ be two local \dr s of the same hypersurface germ
and $E_1,E_2\sset\cOn^{\oplus d}$ the corresponding embedded kernel modules.
Then
\beq \cM_1=\cM_2 \text{ or } \cM_1=A\cM_2 \text{ or } \cM_1=A\cM_2B, \text{ for }A,B \text{ locally invertible on }(\k^n,0)
\eeq
iff
\beq\label{Eq.Notions.of.Equality.of.Kernels}
\Big(E_1,\{s^1_1..s^1_d\}\Big)=\Big(E_2,\{s^2_1..s^2_d\}\Big)\sset\cOn^{\oplus d}  \text{ or }
E_1=E_2\sset\cOn^{\oplus d}  \text{ or } E_1\approx E_2
\eeq
\\2. In particular, if two kernel modules of the same hypersurface are abstractly isomorphic then
their isomorphism is induced by a unique ambient automorphism of $\cOn^{\oplus d}$.
\\3. $\cM$ is decomposable (or equivalent to an upper-block-triangular form) iff $E$
is a direct sum (or an extension).
\\4. Let $\cM_1,\cM_2\in Mat(d\times d,|\cO_{\P^n}(1)|)$ be two global \dr s of the same hypersurface, for $n>1$.
 Let $E_1,E_2$ be the corresponding kernel sheaves. Then the global versions of all the statements above hold.
\eprop
\beR
$\bullet$ In equation (\ref{Eq.Notions.of.Equality.of.Kernels}), in the first case the coincidence of
the natural bases is meant, in the second case the
coincidence of the embedded modules, in the third case the abstract isomorphism of modules.
\li
Part 2 of the last proposition does not hold for arbitrary modules (not kernels). For example
the ideals $<x^l>\sset\k[x]$ for $l\ge0$, are all abstractly isomorphic as (non-embedded) modules
but certainly not as ideals, i.e. embedded modules.
\li
Note that the coincidence/isomorphism of kernel sheaves is a much stronger property than the pointwise
coincidence of kernels as embedded vector spaces. For example, let $\cM$ be a local \dr\  of the
plane curve $C=\{f(x,y)=0\}$. Let $v_1,..,v_p$ be the columns of $\cMv$, i.e. the generators of the
kernel $E_C$. Let $\{g_i=0\}_{i=1}^p$ be some local curves intersecting $C$ at the origin only.
Then $Span(g_1 v_1,..,g_p v_p)$ coincide pointwise with $Span(v_1,..,v_p)$ as a collection of embedded
vector spaces on $C$. Though the two modules correspond to \dr s of distinct curves.
\eeR
\bpr (of proposition \ref{Thm.Notions.of.equality of kernels})
\\{\bf 0.} For modules. As $(\cM\cMv)|_{(X,0)}=\zero$, the  columns of $\cMv|_{(X,0)}$ generate a submodule of $E_{(X,0)}$.
 For any element $s\in E$, one has $\cM s=\det(\cM)v$, where $v$ is some d-tuple.
 Then $\cM(s-\cMv v)=0$ on $(\k^n,0)$. And $\cM$ is non-degenerate on $(\k^n,0)$, so $s=\cMv v$.

For sheaves: the columns of $\cMv$ generate a subsheaf of $E$. If $s\in H^0(E_X)$ then $s$ is the column
whose entries are sections of $\cO_X(d-1)$. By the surjection $H^0(\cO_{\P^n})\to H^0(\cO_X)\to0$ the entries of $s$
are restrictions of some sections of $\cO_{\P^n}$. Hence $s$ is the restriction of some globally defined section $S$,
for which: $\cM S=\det(\cM)(..)$. But then $S=\cMv(..)$. Hence $s$ belongs to the span of the columns of $\cMv$,
thus $h^0(E_X)=d$, i.e. $H^0(E_X)$ is generated by the columns of $\cMv$.

{\bf 1, 2.} As the kernel is spanned by the columns of $\cMv$ the statement is straightforward,
except possibly for the last part: if $E_1\approx E_2$ then $\cM_1=A\cM_2B$.

Let $\phi: E_1\isom E_2$ be an abstract isomorphism, i.e. an $\cO_{(X,0)}$-linear map.
This provides an additional minimal free resolution of $E_1$:
\beq\bM
0&\to &E_1&\to& \cO^{\oplus d}_{(X,0)}&\stackrel{\cM_1}{\to}&\cO^{\oplus d}_{(X,0)}...\\
&&\phi\downarrow&&\psi\downarrow&&\\
0&\to &E_2&\to& \cO^{\oplus d}_{(X,0)}&\stackrel{\cM_2}{\to}&\cO^{\oplus d}_{(X,0)}...
\eM\eeq
By the uniqueness of minimal free resolution, \cite[\S20]{Eisenbud-book}, we get that $\psi$ is an isomorphism.

{\bf 3.} Suppose $E=E_1\oplus E_2$, let $F_2\stackrel{\cM}{\to} F_1\to E\to 0$ be the minimal resolution.
Let $F^{(\al)}_2\stackrel{\cM_\al}{\to} F^{(\al)}_1\to E_\al\to 0$ be the minimal resolutions of $E_1,E_2$. Consider their direct sum:
\beq
F^{(1)}_2\oplus F^{(2)}_2\stackrel{\cM_1\oplus\cM_2}{\to} F^{(1)}_1\oplus F^{(2)}_1\to E_1\oplus E_2=E\to 0
\eeq
This resolution of $E$ is minimal. Indeed, by the decomposability assumption the number of generators of $E$ is the sum of
those of $E_1,E_2$, hence $rank(F_1)=rank(F^{(2)}_1)+rank(F^{(1)}_1)$.  Similarly, any linear relation between the generators of $E$
(i.e. a syzygy)
 is the sum of relations for $E_1$ and $E_2$. Hence $rank(F_2)=rank(F^{(2)}_2)+rank(F^{(1)}_2)$.

Finally, by the uniqueness of the minimal resolution we get that the two proposed resolutions of $E$ are isomorphic, hence the statement.

{\bf 4.}
The statement, $E_1\approx E_2$ implies $\cM_1=A\cM_2B$, is proved for sheaves in \cite[theorem 1.1]{Cook-Thomas79}.
 Note that in general it  fails for $n=1$, see \cite[pg.425]{Cook-Thomas79}.

From this part the rest of the statements follow.
\epr

\subsubsection{\Mg\dr s}\label{Sec.Preliminaries.Maximaly.Generated.DR.s}
Note that at each point $corank\cM |_{pt}\le mult(X,pt)$ (see property \ref{Thm Localiz Chip off Unity}).
This motivates the following
\bed
$\bullet$ A \dr\ of a hypersurface is called \underline{\mg} at the point $pt\in X$ (or Ulrich-maximal \cite{Ulrich84})
 if $corank\cM |_{pt}=mult(X,pt)$.
\li A \dr\ of a hypersurface is called \mg\ \underline{near} the point if it is \mg\ at each point of some
neighborhood of $pt\in X$.
\li A \dr\ of a hypersurface is called \underline{generically} \mg\ if  it is \mg\ at the generic smooth point
of $X_{red}$.
\eed
\bex\label{Ex.Max.Generated.Det.Reps}
1. The \dr s in example \ref{Ex.Det.Reps.of.A_k} and in the first part of example \ref{Ex.Kernel.Maps.and.Images}
are \mg. In fact, as follows from property \ref{Thm.Max.Gen.Are.Decomposable} below, for the ordinary multiple point
 (i.e. curve singularity with several smooth pairwise non-tangent branches) the diagonal matrix is the
{\em unique} local \mg\ representation.
\\2. Any \dr\ of a smooth hypersurface is \mg\ and any \dr\ of a reduced hypersurface is generically \mg.
If the (reduced) hypersurface germ has an isolated singularity then any \dr\ is \mg\ on the punctured neighborhood
of the singular point.
\\3. If $\cM$ is \mg\ at the generic smooth point of $X_{(red)}$ then it is \mg\ at any smooth of $X_{(red)}$,
as the corank of the matrix does not increase under deformations.
\eex
\Mg\ \dr s  are studied in \cite{Kerner-Vinnikov2010}.
They possess various excellent properties, in particular tend to be decomposable:
\bpro\label{Thm.Max.Gen.Are.Decomposable}\label{Thm.Max.Generated.Is.Divisible}
{\bf 1.} Let $\cM$ be a \dr\ \mg\ at the generic smooth points of the (non-reduced)
hypersurface $\{\prod f^{p_\al}_\al=0\}\sset(\k^n,0)$.
Then any entry of $\cMv$ is divisible by $\prod f^{p_\al-1}_\al$.
\\
{\bf 2.} For the case $n=2$, plane curves. Let $(X,0)=(X_1,0)\cup(X_2,0)\sset(\k^2,0)$, here $\{(X_\al,0)\}_\al$ can
be further reducible, non-reduced but without common components.
 Let $\cM$ be a \mg\ local \dr\ of $(X,0)$. Then $\cM$
is locally equivalent to $\bpm \cM_1&*\\\zero&\cM_2\epm$, where $\cM_\al$ are \mg\ \dr s of $(X_\al,0)$.

If in addition the curve germs $(X_i,0)$ have no common tangents then any \mg\ \dr\ is decomposable,
$\cM=\cM_1\oplus\cM_2$.
\\
{\bf 3.} For the case $n>2$. Let $(X,0)=(X_1,0)\cup(X_2,0)\sset(\k^n,0)$, let $\cM_{d\times d}$ be a \dr\ of $(X,0)$
and $E$ its kernel module. Let $E|_{(X_i,0)}/Torsion$ be the restriction to a component. Suppose it
is minimally generated by $d_i$ elements. Similarly,suppose $E^l|_{(X_i,0)}/Torsion$ is minimally generated
by $d^l_i$ elements.

The \dr\ is equivalent to an upper block-triangular iff $d_1+d^l_2=d$ or $d^l_1+d_2=d$.
\epro
\subsubsection{\XXS \dr s}\label{Sec.Preliminaries.XXS.DR.s}
Let $X'\norm X$ be a finite modification, i.e. $X'$ is a pure dimensional scheme, $\nu$ is a
finite surjective proper morphism that is an isomorphism over $X\smin Sing(X_{red})$. If $X$ is reduced
then the "maximal" modification is the normalization $\tX\to X$ and all other modifications
are "intermediate", $\tX\to X'\to X$.
\bed\label{Def.Saturated.Det.Reps}
A \dr\ $\cM$ of a hypersurface is called \underline{\XXS}   if
every entry of $\cMv$ belongs to the relative adjoint ideal
\[
Adj_{X'/X}=\{g\in\cO_X| \nu^*(g)\cO_{X'}\sset\nu^{-1}\cO_X\}
\]
\eed
Any \dr\ is \XXS for the identity morphism $X'\isom X$. A \dr\ is \XXS iff its kernel is a \XXS module, i.e.
$E=\nu_*(\nu^*E/Torsion)$.

As in the \mg\ case, the  \XXS \dr s possess various excellent properties, in particular tend to be decomposable.
\bpro\label{Thm.X'/X-sat.are.decomposable}\cite{Kerner-Vinnikov2010}.
1. Consider the modification
\beq
(X',0)=(X_1,0)\coprod(X_2,0)\norm(X,0)=(X_1,0)\cup(X_2,0)
\eeq
which is the separation of the components. Then $\cM$  is decomposable, $\cM\sim\cM_1\oplus\cM_2$, iff it is \XXS.
\\2. In particular if $(X,0)=\cup_\al(X_\al,0)$ is the union of smooth hypersurface germs and
$(X',0)=\coprod_\al(X_\al,0)$, then any \XXS \dr\ of $(X,0)$ is equivalent to the diagonal matrix (in particular it is \mg).
\\3. If the \dr\ $\cM$ of a plane curve is $\tC/C$-saturated, where $\tC\norm C$ is the normalization,
then $\cM$ is \mg\ and the entries of $\cMv$ generate the adjoint ideal $Adj_{\tC/C}$.
\epro
Finally we state an additional decomposability criterion from \cite{Kerner-Vinnikov2010}:
\bpro\label{Thm.Kernel.Has.Independent.Fibres.Decomposable.}
Let $(X,0)=\cup(X_\al,0)\sset(\k^n,0)$ be a collection of reduced, smooth hypersurfaces.
The \dr\ $\cM$ of $(X,0)$ is completely decomposable iff the geometric fibres $\{E_\al|_0\}$
are linearly independent: $Span(\cup E_\al|_0)=\oplus E_\al|_0$.
\epro
\section{Global decomposability}\label{Sec Global Decomposability}
The local decomposability at each point implies the global one.
\bthe\label{Thm Decomposability Global from Decomposability Local}
Let $X=X_1\cup X_2\sset\P^n$ be a global decomposition of the hypersurface. Here $X_1,X_2$ can
be further reducible, non-reduced, but without common components.
 Then $\cM$ is globally decomposable, i.e. $\cM\stackrel{globally}{\sim}\cM_1\oplus\cM_2$,
 iff it is locally decomposable at each point $pt\in X$, i.e.
$\cM\stackrel{locally}{\sim}\one\oplus\cM_1|_{(\P^n,pt)}\oplus\cM_2|_{(\P^n,pt)}$.
Here $\cM_\al|_{(\P^n,pt)}$ are the local \dr s near $pt\in\P^n$, one works over $\cOn=\k[x_1,..,x_n]_{(\cm)}$.
\ethe
\bpr $\Lleftarrow$ Let $f=f_1f_2$ be the homogeneous polynomials defining $X,X_1,X_2$. Here $f_\al$ can be reducible,
non-reduced,
but mutually prime.
\\{\bf Part1.} By the assumption at each point $\cMv\stackrel{locally}{\sim}f\one\oplus f_2\cMv_1\oplus f_1\cMv_2$.
The local ideal of $\cOn$  generated by the entries of $\cMv$ is invariant under local equivalence. Hence we get:
any entry of $\cMv$ at any point $pt\in\P^n$ belongs to the
local ideal $\bl f_1,f_2\br\sset\cOn$.

Now use the Noether's $AF+BG$ theorem \cite[pg. 139]{ACGH-book}:\\
{\em Given some homogeneous polynomials $F_1..F_k$,
 whose zeros define a subscheme of $\P^n$ of dimension $(n-k)$.
Suppose at each point of $\P^n$ the homogeneous polynomial $G$ belongs to the local ideal generated
by $F_1..F_k$. Then $G$ belongs to the global ideal in $\k[x_0,..,x_n]$ generated by $F_1..F_k$.}

Apply this to each entry of $\cMv$. Hence we get, in the matrix notation:
\beq
\cMv=f_2\cNv_1+ f_1\cNv_2, \hspace{1cm} \cNv_\al\in Mat\Bigg(d\times d,H^0\big(\cO_{\P^n}(deg(f_\al)-1\big)\Bigg)
\eeq
\\
{\bf Part2.} From the last equation we get: $f_1f_2\one=\cM\cMv=f_2\cM \cNv_1+f_1\cM \cNv_2$.
 Note that $f_1,f_2$ are relatively prime, thus:
$\cM\cNv_\al=f_\al A_\al$. Here $A_\al$ is a $d\times d$ matrix whose entries are forms of degree zero,
i.e. {\em constants}. Similarly, by considering $\cMv\cM$ we get $\cNv_\al\cM=f_\al B_\al$.
\\\\
Note that $A_1+A_2=\one=B_1+B_2$, by their definition. In addition $A_\al A_\be=\zero=B_\al B_\be$ for $\al\neq \be$.
Indeed: $f_\al B_\al \cNv_\be=\cNv_\al \cM \cNv_\be=f_\be\cNv_\al A_\be$. So $B_\al \cNv_\be$ is divisible by $f_\be$
and
 $\cNv_\al A_\be$ is divisible by $f_\al$. But $\{f_\al\}$ are mutually prime and the degree
 of the entries in $\cNv_\al $ is $d_\al-1$.
Therefore: $B_\al \cNv_\be=0=\cNv_\al A_\be$ (for $\al\neq\be$). And this causes $B_\al (\cNv_\be\cM )=0=(\cM \cNv_\al A_\be)$.
\\\\
Thus $\{A_\al \}$ and $\{B_\al \}$ form a partition of identity, i.e. $\oplus A_\al =\one=\oplus B_\al $.
So, one  can bring the collections $\{A_\al \}$, $\{B_\al \}$ to the block-diagonal form:
\beq
\bpm \tA_1&0\\0&\tA_2\epm=\one=\bpm \tB_1&0\\0&\tB_2\epm
\eeq
This is done by the multiplication $\cM \to U_1\cM U_2$ (and accordingly $\cMv\to U^{-1}_2\cMv U^{-1}_1$),
which acts on $A,B$ as: $A_\al \sim\cMv_\al \cM \to U^{-1}_2A_\al U_2$ and $B_\al \sim\cM \cMv_\al \to U_1B_\al U^{-1}_1$.

Then from the definitions $\cNv_\al \cM \sim B_\al $ and $\cM \cNv_\al \sim A_\al $ one gets: $\cNv_\al \sim B_\al \cMv$ and
$\cNv_\al \sim\cMv A_\al $. So, $\cNv_\al $ is just one block on the diagonal. Thus $\cMv$ is block diagonal.

Finally note that from $\cNv_\al \cM =f_\al  A_\al$ it follows that $\det(\cNv_\al )f_\al =f^{d_\al }_\al \times$const.
So the multiplicities are determined uniquely.
\epr
If we combine the theorem with property \ref{Thm.Max.Gen.Are.Decomposable},
and property \ref{Thm.X'/X-sat.are.decomposable} we get:
\bcor
Let $X=\cup p_\al X_\al \sset\P^n$ be the global decomposition into distinct irreducible reduced projective hypersurfaces.
\\1. Let $X'=\coprod p_\al X_\al \norm \cup p_\al X_\al $ be the separation of components (a finite modification).
Any \XXS \dr\  is globally completely decomposable, $\cM\sim\oplus\cM_\al $, where $\cM_\al $ is a \dr\ of $p_\al X_\al $.
\\2. Suppose for each point $pt\in X$ any two components passing through this point $pt\in X_\al \cap X_\be$
intersect generically transverse, i.e. $(X_\al \cap X_\be)$ is reduced. Then any \mg\ \dr\ of $X$
is globally completely decomposable, $\cM\sim\oplus\cM_\al $.
\ecor

It is interesting that this local-to-global theorem, a non-trivial result in linear algebra,
is immediate if viewed as a statement about the kernel sheaves.
\bprop\label{Thm Decomposability Global From Sheaf Axioms}
 Let $X=X_1\cup X_2\sset\P^n$, here $X_\al $ can be further reducible, non-reduced, but with no common components.
Let $\cM$ be a global \dr\  of $X$.
\\{\bf 1.} Suppose at each point $pt\in X_1\cap X_2$ the \dr\  is locally decomposable,
i.e.  $\cM\stackrel{loc}{\sim}\cM_1\oplus\cM_2$ with $\cM_\al $ the local \dr\  of $X_\al $.
Then $\cM$ is globally decomposable.
\\{\bf 2.} Suppose at each point $pt\in X_1\cap X_2$ the \dr\  can be brought locally to the upper block-triangular form,
i.e.  $\cM\stackrel{loc}{\sim}\bpm\cM_1&..\\\zero&\cM_2\epm$ with $\cM_\al $ the local \dr\  of $X_\al $.
Then $\cM$ is globally equivalent to an upper-block-triangular matrix.
\eprop
\bpr For $n=1$ any \dr\  is completely decomposable, hence we assume $n>1$ and $X$ is connected.
\\{\bf 1.} Let $E$ be the kernel sheaf of $\cM_X$, let $X_\al\stackrel{i_\al}{\hookrightarrow}X$ be the natural embeddings.
Define the restriction to the component: $E_\al:=i^*_\al(E)/Torsion$.

Consider the new sheaf on $X$: $E':=(i_1)_*(E_1)\oplus (i_2)_*(E_2)$. It is a coherent sheaf and by construction
there is a natural globally defined map: $E\to E'$, the direct sum of two restrictions. We claim that
it is isomorphism on all the stalks. It is an isomorphism on $X\smin(X_1\cap X_2)$, so the kernel of this map would
be a torsion subsheaf of $E$. But $E$ is torsion-free, hence the map is injective. So, one only need to check
the surjectivity (over $X_1\cap X_2$), which holds by construction.

Therefore by the axiom of sheaf the global map is an isomorphism too:
\beq
E\isom E'= i^*_1(E_1)\oplus i^*_2(E_2)
\eeq
So, by the global automorphism of the sheaf,
i.e. a linear operator with constant coefficient, $\cM$ is brought to the block diagonal form, see
proposition \ref{Thm.Notions.of.equality of kernels}.

{\bf 2.} Define the sheaf $E'$ on $X$ as the extension
\beq
0\to (i_1)_*(E_1)\to E'\to (i_2)_*(E_2)\to0
\eeq
 where at each intersection point of $X_1\cap X_2$ the stalk of $E'$ is defined by its extension class:
$[E']=[E]\in Ext\big((i_2)_*(E_2),(i_1)_*(E_1)\big)$. Again one has the natural map $E\to E'$, which is an
isomorphism on stalks, hence a global isomorphism.
\epr
\bex
Consider a reducible (reduced) curve $X=X_1\cup X_2\sset\P^2$, suppose all the intersections of $X_1,X_2$ are transverse,
i.e. the corresponding singularities are nodes. Let $\cM$ be a \dr\ of $X$. Then it either splits, $\cM\sim\cM_1\oplus\cM_2$,
or for at least one intersection point, $pt\in X_1\cap X_2$, the corank is minimal: $corank(\cM|_{pt})=1$.
\eex
\section{(Co-)Kernel sheaves of \dr s}\label{Sec.Kernels.Properties.Classification}
\subsection{Properties and classification of kernels on the hypersurface}
Recall (from \S\ref{Sec.Preliminaries.Kernel.Sheaves}) that the (co)kernel sheaves are defined by:
\beq\label{Eq Exact Sequences For Kernel Sheaves. Section 4}
\ber
0\to E_X \to\cO_X^{\oplus d}(d-1)\stackrel{\cM}{\to}\cO_X^{\oplus d}(d)\to Coker(M)_X\to0,\\
0\to E_X^l \to\cO_X^{\oplus d}(d-1)\stackrel{\cM^T}{\to}\cO_X^{\oplus d}(d)\to Coker(M^T)_X\to0
\eer\eeq
\bthe\label{Thm.Kernel.Sheaves.Properties}
{\bf 1.} The sheaf $E_X$ has periodic resolutions:
\[
...\stackrel{\cM}{\to}\cO_X^{\oplus d}(-d)\stackrel{\cMv}{\to}\cO_X^{\oplus d}(-1)\stackrel{\cM}{\to}
\cO_X^{\oplus d}\stackrel{\cMv}{\to}E_X\to 0,\]
\[0\to E_X\to\cO_X^{\oplus d}(d-1)\stackrel{\cM}{\to}\cO_X^{\oplus d}(d)\stackrel{\cMv}{\to}\cO_X^{\oplus d}(2d-1)..\]
In particular, for any point $pt\in X$, the stalk $E_{(X,pt)}$ is a Cohen-Macaulay module over $\cO_{(X,pt)}$.
\\{\bf 2.}  $h^0(E_X(-1))=0$. $h^i(E_X(j))=0$, for $0<i<n-1$ and $j\in\Z$.  $h^{n-1}(E_X(-j))=0$, for $j<n$.
In particular $\chi(E(-j))=0$ for $j=1,..,n-1$.
\\{\bf 3.} The sheaves $E_X$ and $Coker(M)_X$ are torsion free.
{\\\bf 4.} The sheaf $E_X$ is quasi-locally free iff $dim(Ker\cM|_{pt})=const$ on $X$.
The sheaf $E_X$ is locally free iff $dim(Ker\cM|_{pt})=1=const$ on $X$.
In particular if the hypersurface is smooth the kernel sheaf is invertible.
{\\\bf 5.} Suppose $X\sset\P^n$ is reduced and $E_X$ is locally free. Then $n<4$ and
\[c_1(E_X)L^{n-2}=\frac{(d-1)}{2}L^{n-1},\quad c^2_1(E_X)L^{n-3}=\frac{(d-1)(d-2)}{6}L^{n-1}\]
where $[L]\in H^2(X,\Z)$ is the hyperplane class and $\deg(L^{n-1})=d$, the degree of the zero-dimensional cycle.
\\\\
All the properties hold for $E^l_X$  (is generated by the columns of  $(\cMv|_X)^T$ etc).
\ethe
\bpr {\bf 1.} Follows immediately from property \ref{Thm.Notions.of.equality of kernels}.

{\bf 2.} Let $I_X\sset\cO_{\P^n}$ be the defining ideal of the hypersurface, i.e. $0\to I_X\to\cO_{\P^n}\to\cO_X\to0$.
Note that $I_X\approx\cO_{\P^n}(-d)$ and $h^{i>0}(\cO_{\P^n}(j))=0$ for $i\neq n$ or $j+n+1>0$.
Therefore the long exact sequence of cohomology gives:
\beq
h^i(\cO_X(j))=0\quad\text{for }0<i<n-1,\ j\in\Z,\quad\text{and }h^{n-1}(\cO_X(j))=0\quad\text{for }j>d-n-1
\eeq
In addition $H^0(\cO_{\P^n}(j))\to H^0(\cO_X(j))\to0$ for $n>2$ or for $n=2$, $j>d-3$.
\\\

Now we prove that $h^0(E_X(-1))=0$. Let $0\neq s\in H^0(E_X(-1))$,
i.e. $s$ is a d-tuple whose entries are sections of $\cO_X(d-2)$. As
 $H^0(\cO_{\P^n}(d-2))\to H^0(\cO_X(d-2))\to0$, we get a global section $S\in H^0(\cO^{\oplus d}_{\P^n}(d-2))$
that restricts to $s$ on the hypersurface. Thus $\cM_{\P^n}S$ is proportional to $f$.
But the total degree of any entry of $\cM_{\P^n}S$ is $(d-1)$, which is less than $deg(f)=d$. So $\cM_{\P^n}S=0$,
contradicting the global non-degeneracy of $\cM$.
\\\

The vanishing $h^{i>0}(E_X(j))=0$. Use the exact sequence
\beq
0\to E_X(-1)\to \cO_X^{\oplus d}(d-2)\stackrel{\cM}{\to}Im(\cM)\to0
\eeq
Twist it by $\cO(-j)$, for $j\ge0$, note that $h^0(E_X(-1-j))=0$. Further, note that the map
$H^0(\cO_X^{\oplus d}(d-2-j))\to H^0(Im(\cM)(-j))$ is surjective. Indeed, $Im(\cM)\sset\cO_X^{\oplus d}(d-1)$,
hence any global section of $Im(\cM)(-j)$ is the restriction of some section of $\cO_{\P^n}^{\oplus d}(d-j-1)$.
The sequence of cohomologies consists of the parts:
\beq\label{Eq.Seq.of.Cohoms.1}
H^i(\cO_X^{\oplus d}(d-2-j))\to H^i(Im(\cM)(-j))\to H^{i+1}(E_X(-1-j))\to H^{i+1}(\cO_X^{\oplus d}(d-2-j))
\eeq
For $n=2$ the case $j=0=i$ gives: $h^{1}(E_X(-1))=0$, proving the statement.
\\
\\\

Suppose $n>2$. Then the last sequence for $i=0$ and $j\ge0$ gives: $H^{1}(E_X(-1-j))=0$.
Further, the sheaf $Im(\cM)$ can be considered as the kernel sheaf of
$\cO_X^{\oplus d}(d-1)\stackrel{\cMv}{\to}\cO_X^{\oplus d}(2d-2)$.
Therefore we have the additional exact sequence
\beq
0\to Im(\cM)\to \cO_X^{\oplus d}(d-1)\stackrel{\cMv}{\to}E_X(d-1)\to0
\eeq
Twist it and take cohomology:
\beq\label{Eq.Seq.of.Cohoms.2}
H^i(\cO_X^{\oplus d}(-j))\to H^{i}(E_X(-j))\to H^{i+1}(Im(\cM)(1-d-j))\to  H^{i+1}(\cO_X^{\oplus d}(-j))
\eeq
Similarly to $H^{1}(E_X(-1-j))=0$ one gets $H^{1}(Im(\cM)(-j))=0$. These are the ''initial conditions''.
Combine now equations (\ref{Eq.Seq.of.Cohoms.1}) and (\ref{Eq.Seq.of.Cohoms.2}) with the initial conditions to get:
\beq
H^{i}(E_X(-j))=0,\quad\quad H^{i}(Im(\cM)(-j))=0,\quad\text{for }0<i<n-1,\quad j\ge0
\eeq
while for $i=n-1$ one has:
\beq
0\to H^{n-1}(E_X(-1-j))\to  H^{n-1}(\cO_X^{\oplus d}(d-2-j))\to H^{n-1}(Im(\cM)(-j))\to 0
\eeq
giving $H^{n-1}(E_X(-1-j))=0$ for $j<n-1$.
\\{\bf 3.} The sheaves $E,E^l$ are torsion-free as sub-sheaves of the torsion-free sheaf $\cO_X^{\oplus d}(d-1)$.

To check that $Coker(\cM|_X)$, $Coker(\cM^T|_X)$ are torsion free it is enough to consider their stalks at a point.
Then for $s\in \cO_{(X,0)}^{\oplus d}/\cM\cO_{(X,0)}^{\oplus d}$ and $g\in\cO_{(X,0)}$ not a zero divisor we
should prove:
if $gs=0\in\cO_{(X,0)}^{\oplus d}/\cM\cO_{(X,0)}^{\oplus d}$ then $s\in\cM\cO_{(X,0)}^{\oplus d}$.

But if $gs\in\cM\cO_{(X,0)}^{\oplus d}$ then
$\cMv gs\in\cMv\cM\cO_{(X,0)}^{\oplus d}=\det(\cM)\cO_{(X,0)}^{\oplus d}\equiv 0$. So $\cMv_X s=0$, hence $s=\cM s_1$.
\\{\bf 4.}
Suppose $E_X$ is quasi-locally free at $0\in(X,0)$, i.e. the stalk of its associated graded sheaf
(\S\ref{Sec.Preliminaries.Sheaves.on.Sing Curves}) is a free $\cO_{(X_{red},0)}$ module.
We can assume $\cM|_0=\zero$. Let $s_1,..,s_k$ be the generators of $E_{(X,0)}$,
hence if $\sum g_i s_i=0$ for some $\{g_i\in \cO_{(X_{red},0)}\}$ then $\{g_i=0\in\cO_{(X_{red},0)}\}$.

But $\{s_1,..,s_k\}$ are columns of $\cMv$, so $(s_1,.,,s_k)\cM=\det(\cM)\one=\zero$. Hence
$\cM|_{(X_{red},0)}=\zero$. In particular $corank(\cM|_X)=const$.

If $E_X$ is locally free then by the same argument  get $\cM|_{(X,0)}=\zero$,
 which is possible only if $\cM$ is a $1\times 1$ matrix.
\\{\bf 5.} If the kernel sheaf is locally free then its fibre at each point of $X$ is of dimension one. Hence at
the points where $corank\cM>1$ the kernel cannot be locally free. We claim that the locus of such point is of
codimension at most 4 in $\P^n$. This locus is defined by the vanishing of all the maximal minors of $\cM$. Hence
we can use the standard fact: inside the parameter space of all the $d\times d$ matrices (with constant coefficients)
the locus of matrices of corank at least two has codimension four.

Assume local freeness of $E_X$, and hence $n<4$. As was proved above $\chi(E(-j))=0$ for $j=1,..,n-1$.
Apply Hirzebruch-Riemann-Roch theorem, \S\ref{Sec.Preliminaries.Sheaves.on.Sing Curves}. For $n=2$ we get:
\beq
0=\chi(E(-1))=deg(E(-1))-(p_a-1),\quad\rightsquigarrow\quad deg(E)=\frac{d(d-1)}{2}
\eeq
For $n=3$ we get:  $\Big(ch(E(-1-j))Td(X)\Big)_{top.dim.}=0$ for $j=0,1$.
Note that $ch(E(-1-j))=ch(E)ch\big((-1-j)L\big)$. Hence we get the system of two equations:
\beq
\Big((\sum_{i\ge0}\frac{((-1-j)L)^i}{i!})ch(E)Td(X)\Big)_{top.dim.}=0,\quad j=0,1.
\eeq
The first bracket produces the matrix: $\bpm 1& -L&\frac{(-L)^2}{2!}\\1&-2L&\frac{(-2L)^2}{2!}\epm$.
Apply row operations to bring this matrix to the form:
$\bpm 1& 0&-L^2\\0&-L&\frac{3L^2}{2}\epm$.
Now substitute this to the initial equations to get:
\beq
L(c_1(E_X)+\frac{c_1(T_X)}{2})=\frac{3L^2}{2},\quad
(\frac{c^2_1(E_X)}{2}+c_1(E_X)\frac{c_1(T_X)}{2}+\frac{c^2_1(T_X)+c_2(T_X)}{12})=L^2
\eeq
Finally, use example \ref{Ex.Chern.Class.Hypersurface}.
\epr
\beR
The vector bundles satisfying $h^{0<i<dim(X)}(E_X(j))=0$ are called ACM (arithmetically Cohen-Macaulay) bundles.
For example on smooth quadrics in $\P^n$ there are only two such bundles (up to a twist): the trivial and the
spinor bundles. For a recent review cf. \cite{Pons-Llopis.Tonini2009}.
\eeR
It is possible to give a simple criterion for a torsion free sheaf to be the kernel of some \dr.
\bthe\label{Thm.Kernel.Sheaves.Classification}
Let $X=\cup p_\al X_\al\sset\P^n$ be a hypersurface of degree $d$. Let $E_X$ be a torsion-free sheaf
of multi-rank $\{p_\al\}_\al$, suppose at each point of $X$ the stalk $E|_{(X,pt)}$ is a maximally
Cohen-Macaulay module over $\cO_{(X,pt)}$.
Assume $h^0(E_X(-1))=0$, $h^i(E_X(j))=0$ for $0<i<n-1$ and $j\in\Z$, and $h^{n-1}(E_X(1-n))=0$.
 Then $E_X$ is the kernel sheaf of a \dr\ of $X$, i.e.
\beq
0\to E_X \to\cO_X^{\oplus d}(d-1)\stackrel{\cM}{\to}\cO_X^{\oplus d}(d)\to\dots
\eeq
where $\cM$ is a \dr\ of $X\sset\P^n$.
\ethe
\bpr
{\bf Step 1.} The preparation.
 Define the auxiliary sheaf: $E^l_X:=E_X^*\otimes w_X(n)$. As $X$ is Gorenstein, the dualizing sheaf
 $w_X$ is invertible, hence $E^l_X$ is torsion-free.
By Serre duality: $h^i(E^l_X(j))=h^{n-1-i}(E_X(-n-j))$. Hence $h^i(E^l_X(j))=0$ for $0<i<n-1$, $j\in\Z$ and
$h^0(E^l_X(-1))=h^{n-1}(E_X(1-n))=0$ and $h^{n-1}(E^l_X(1-n))=h^0(E_X(-1))=0$.

Now we prove that $h^0(E^l_X)=d=h^0(E_X)$. If $n=2$, i.e. $X$ is a curve, and $E_X$ is invertible, then this
follows from Riemann-Roch:
\beq
\chi(E_X)-\chi(E_X(-1))=\deg(E_X)-\deg(E_X(-1))=\deg(\cO_X(1))=d
\eeq
In higher dimensions one can argue as follows.
Let $L\sset\P^n$ be a line. We assume $L$ is generic enough, such that it intersects the reduced hypersurface $X_{red}$
at  smooth points only and the stalks of $E$ at $Z=L\cap X$ are quasi-locally free
(cf.\S\ref{Sec.Preliminaries.Sheaves.on.Sing Curves}). Let $I_Z\sset\cO_X$
be the defining sheaf of ideals:
\beq
0\to I_Z\to\cO_X\to\cO_Z\to0
\eeq
The line $L$ is the intersection of $(n-1)$ hyperplanes, hence $I_Z$ admits Koszul resolution.
Tensor the exact sequence above with $E_X$:
\beq
Tor^1_{\cO_X}(I_Z,E_X)\to \underbrace{Tor^1_{\cO_X}(\cO_X,E_X)}_{=0}\to Tor^1_{\cO_X}(\cO_Z,E_X)\to E_X\otimes I_Z\to E_X\to E_X\otimes\cO_Z\to0
\eeq
We claim that $Tor^1_{\cO_X}(\cO_Z,E_X)=0$. Indeed, let $\to P_2\to P_1\to E_X\to0$ be a projective resolution, then
$Tor^1_{\cO_X}(\cO_Z,E_X)$ is the cohomology of the tensored complex at the place $\to P_2\otimes\cO_Z\to$.
As $Z$ consists of several generic points of $X$, the latter complex is exact, hence the cohomology vanishing.
Then from the exact sequence above
we get also:  $Tor^1_{\cO_X}(I_Z,E_X)=0$.

Finally, take the cohomology of the sequence $0\to E_X\otimes I_Z\to E_X\to E_X\otimes\cO_Z\to0$.
We claim that $h^1(E_X\otimes I_Z)=0$. Indeed, the resolution of $I_Z$ is: $0\to Syz\to\cO_X^{\oplus d}(-1)\to I_Z\to0$,
where $Syz$ is the syzygy module. Multiply this by $E_X$, as mentioned above $Tor^1_{\cO_X}(I_Z,E_X)=0$, thus:
$0\to Syz\otimes E_X\to\cO_X^{\oplus d}(-1)\otimes E_X\to I_Z\otimes E_X\to0$.
Now take the cohomology, use $h^i(E_X(-1))=0$ to get: $h^i(I_Z\otimes E_X)=h^{i+1}(Syz\otimes E_X)$.
As $I_Z$ has Koszul resolution the resolution of syzygy module is Koszul too (and shorter). Hence we get a
reduction in codimension. Finally, using $h^i(E_X^{\oplus(n-1)}(-1))=(n-1)h^i(E_X(-1))=0$ we have:
$h^1(I_Z\otimes E_X)=0$.

Therefore:
\beq
h^0(E_X)=h^0(E_X\otimes\cO_Z)=length(E_X\otimes\cO_Z)=\suml_{pt\in Supp(Z)}mult(X,pt)=d
\eeq
Here we used proposition \ref{Thm.Rank.of.Sheaf.Equals.Intersection.Multiplicity}: the local rank of $E$ at the
smooth point of $X_\al$ is $p_\al$.

Similarly for $h^0(E^l_X)=d$.
\\
\\\

{\bf Step 2.} Construction of the candidate for $\cMv$.
Consider the global sections $H^0(E_X,X)=Span(s_1..s_{d})$ and $H^0(E^l_X,X)=Span(s^l_1..s^l_{d})$.
By the construction of $E_X,E^l_X$ one has the pairing:
\beq
\Big(H^0(E_X),H^0(E^l_X)\Big)\to H^0(w_X(n))\approx H^0(\cO_X(d-1))
\eeq
The later isomorphism is due to the adjunction for hypersurface in $\P^n$: $w_X\approx\cO_X(d-n-1)$.

By the surjection $H^0(\cO_{\P^n}(d-1))\to H^0(\cO_X(d-1))\to0$, the chosen
basis of global sections defines a matrix of homogeneous polynomials:
\beq\label{Eq Candidate for the adjoint matrix}
<s_i,s^l_j>\to \{A_{ij}\}\in Mat\Big(d\times d,H^0(\cO_{\P^n}(d-1))\Big)
\eeq
The entries $A_{ij}$ are defined up to the global sections of $I_X(d-1)\subset\cO_{\P^n}(d-1)$.
As $I_X\approx\cO_{\P^n}(-d)$ we get $h^0(I_X(d-1))=h^0(\cO_{\P^n}(-1))=0$.
Hence the entries $A_{ij}$ are defined uniquely.
\\
\\\

{\bf Step 3.} We prove that the matrix $A$ is globally non-degenerate.  First consider the case
of reduced hypersurface, so $E,E^l$ are of constant rank 1. It is enough to show that the restriction of $A$ to
a line in $\P^n$ is non-degenerate.

Let $L\sset\P^n$ be the generic line, such that $L\cap X=\{pt_1..pt_d\}$ is the set of distinct reduced points.
Consider the restriction $A|_L$, i.e. a \dr\ of $L\cap X$ in $\P^1$.

By linear transformations applied to $Span(s_1..s_{d})$ and $Span(s^l_1..s^l_{d})$ we can choose
the basis of sections satisfying the conditions:
\beq
L\cap X\supset div(s_i)\ge\sum_{j\neq i}~^{\!\!\! '}pt_j\not\ge pt_i,\quad
L\cap X\supset div(s^l_i)\ge\sum_{j\neq i}~^{\!\!\! '}pt_j\not\ge pt_i
\eeq
Indeed, suppose $s_1|_{pt}\neq0$, then can assume $s_{i>1}|_{pt_1}=0$ and continue with the remaining
points $\{pt_{i>1}\}$ and sections $\{s_{i>1}\}$. At each step there is at least one section that
does not vanish at the given point. Otherwise at the end we get a section vanishing at all the points,
thus producing a section of $E_{L\cap X}(-1)$, contrary to $h^0(E_{L\cap X}(-1))=0$.

By the choice of sections, one has $\forall k$, $\forall i\neq j$: $A_{ij}|_{pt_k}=0$. So, $deg(A_{ij}|_L)\ge d$.
But the entries of $A$ are polynomials on $\P^n$ of degree $(d-1)$. Hence for $i\neq j$: $A_{ij}|_L\equiv0$.
On the contrary: $A_{ii}|_{pt_i}\neq0$, so $A_{ii}|_L\not\equiv0$. Therefore, $A|_{L}$ is a diagonal
matrix, none of whose diagonal entries vanishes identically on the hypersurface $L\cap X$. Thus $\det(A|_L)\not\equiv 0$,
so $\det(A)\not\equiv 0$.
\\
\\
Now consider the non-reduced case: $X=\cup p_\al X_\al$ and the rank of $E_X$ on $X_\al$ is $p_\al$.
The generic line $L$ intersects the hypersurface
along the scheme $\sum_{X_\al\sset X} (p_\al\sum_{j} pt_{\al,j})$. Here $pt_{\al,j}\in L\cap X_\al$ and by
genericity $X_{red}$ is smooth at $pt_{\al,j}$ and
the stalk $E_X$ is quasi-locally free at $pt_{\al,j}$, of rank $p_\al$.

Consider the values of the sections at the point $p_1 pt_{1,1}$, i.e. the $p_1$'th infinitesimal
neighborhood of $pt_{1,1}$.
By linear transformation we can assume: $s_{i>p_1}|_{p_1pt_{1,1}}=0$. In other words these sections
are locally divisible by the defining equation of the germ $(p_1 X_1,pt_{1,1})$.

Similarly for $p_1 pt_{1,2}$: $rank Span\{s_{i>p_1}|_{p_1 pt_{1,1}}\}\le r_1 p_1$ etc.,
going over all the points of the intersections $L\cap X$.

If at least for one point the inequality $rank(Span(..))\le p_\al$ is strict then at the end we get
a non-zero section that vanishes at all the points $p_\al pt_{\al,j}$. But this is a section of $E(-1)$,
contrary to the initial assumptions. Hence at each of the above steps we have equality:
$rank(Span(..))=p_\al$.

Repeating this process we can assume that $s_1,..,s_{p_1}$ vanish at all the points $p_\al pt_{\al, i}$
except for $p_1 pt_{1,1}$. Similarly for all other sections and for the sections of $E^l_X$.
\\
\\
Therefore in the chosen basis, $A|_L$ has a block structure, with diagonal blocks corresponding to the
intersection points $L\cap X$. As in the reduced case one gets that the off diagonal blocks vanish on $L$.
Hence for non-degeneracy of $A|_L$ we need to check the non-degeneracy of all the diagonal blocks.
But each such block has rank $p_\al$ at the point $p_\al pt_{\al,j}$, by the construction above.
So, $A|_L$ is non-degenerate on $L$, hence $A$ is non-degenerate on $\P^n$.
\\
\\{\bf Step 4.} The vanishing orders of the entries of $A|_L$.

Let $p_\al pt_{\al,i}\in L\cap X$. Let $A_{\al,i}$ be the diagonal block corresponding to this point,
i.e. the only block whose entries do not all vanish on $p_\al pt_{\al,i}$.
By the choice of the basic sections as above we have: $(A_{\al,i})_{1,1}$ does not vanish at $pt_{\al,i}$.
$(A_{\al,i})_{1,2}$  and $(A_{\al,i})_{2,1}$ vanish to the first order at $pt_{\al,i}$.
The entries of the next anti-diagonal vanish to the second order etc. Finally the entries below
the main (longest) anti-diagonal, that goes from $(X_{\al,i})_{1,p_\al}$ to $(A_{\al,i})_{p_\al,1}$,
 vanish on $p_\al pt_{\al,i}$.

Hence, any entry of $A^\vee|_L$, i.e. any principal minor of $A|_L$, vanishes at $pt_{\al,i}$ to the order at
least
\beq
(d-1-p_\al)p_\al+p_\al(p_\al-1)=p_\al(d-2)
\eeq
 Here the first summand corresponds to all the blocks
except for $X_{\al,i}$, the second summand corresponds precisely to the main anti-diagonal of the block $A_{\al,i}$.

By going over all the points of $L\cap X$ one has: any entry of $A^\vee|_L$ vanishes on $(d-2)X\cap L$.
As this is true for any generic line $L$ we get: any entry of $A^\vee$ vanishes on $(d-2)X$, i.e. is
divisible by $f^{d-2}$, where $X=\{f=0\}\sset\P^n$.
\\
\\
{\bf Step 5.} Construction of $\cM$.
As any entry of $A^\vee$ is a polynomial of degree $(d-1)^2$ and is divisible by $f^{d-2}$,  we define
$\cM:=\frac{A^\vee}{f^{d-2}}\in Mat(d\times d,H^0(\cO_{\P^n}(1)))$. Then
\beq
\det(\cM)=\frac{\det(A^\vee)}{f^{d(d-2)}}=(\frac{\det(A)}{f^{d-1}})^{d-1}f
\eeq
As the left hand side is a polynomial we get that $det(A)$ is divisible by a power of $\prod f_\al$.
The direct check gives the minimal possible power: $det(A)$ is divisible by  $f^{d-1}$.
But entries of $X$ are of degree $(d-1)$ thus, after scaling by a constant: $\det(A)=f^{d-1}$
and $\det(\cM)=f$.

Finally, by construction $MA=f\one$, while $A=\{(e_i,e^l_j)\}$. Thus $\cM(e_i,e^l_j)\equiv f\zero$,
causing: $\cM e_i\equiv 0(f)$. Therefore $Span(e_1..e_{d})\sset Ker(\cM|_X)$ and $E=Ker(\cM|_X)$
(by the equality of dimensions).
\epr
\subsection{Kernels on the reduced curves/hypersurfaces}\label{Sec.Kernel.Sheaf.On.Reduced.Locus}
Suppose the hypersurface is non-reduced: $X=\cup p_\al X_\al\sset\P^n$, where each $X_\al$ is irreducible reduced.
Then the kernel sheaf on $X$ can be restricted to the sheaf $E_X\underset{\cO_X}{\otimes}\cO_{X_{red}}$ on $X_{red}=\cup X_\al$,
corresponding to the natural embedding of closed points $X_{red}\into X$. In other words one considers
$E$ as the sheaf of modules over $\cO_{X_{red}}$. This is a "non-embedded" restriction.

One could consider an "{\em embedded}" restriction, induced by the projection
$\cO_{X}^{\oplus d}(d-1)\to\cO_{X_{red}}^{\oplus d}(d-1)$.
Namely, take any section of $E_X\sset\cO_{X}^{\oplus d}(d-1)$
and consider its image in $\cO_{X_{red}}^{\oplus d}(d-1)$. Then consider the
sheaf generated by such sections. This naive restriction is far from being injective.
For example, any kernel sheaf of a generically \mg\ \dr\ will go
to zero (as the sheaf is spanned by the columns of $\cMv$ and each of them is divisible by $\prod f^{p_\al-1}_\al$,
property \ref{Thm.Max.Generated.Is.Divisible}).
However, for generically \mg\ \dr s one can define a natural injective reduction.
\bdp
Let $E_X$ be the kernel sheaf of the generically \mg\ \dr\  $\cM$.
Let  $E_X^{red}$ be the sheaf of $\cO_{X_{red}}$ modules, generated by the
 columns of $\frac{1}{\prod_\al f^{p_\al-1}}\cMv$.
 Then $E_X^{red}$ depends on $E_X$ only and has the multi-rank $(p_1,..,p_k)$ on $\cup  X_\al$.
 The reduction $E_X\to E_X^{red}$ is injective on generically \mg\ \dr s.
 The pair $(\cM,\frac{1}{\prod_\al f^{p_\al-1}}\cMv)$ is a matrix factorization of $\prod f_\al$.
\edp
\bpr
By the remark above $\frac{1}{\prod_\al f^{p_\al-1}}\cMv$ is a matrix of polynomials. Further, any equivalence of $\cM$
 induces that of $\cMv$, preserving the isomorphism class of $E_X^{red}$.

Finally, given $E_X^{red}$, generated by the columns of $\cNv$, the kernel $E_X$ must be generated by the columns of
$\prod_\al f^{p_\al-1}\cNv$, i.e. is determined uniquely.
\epr

Similarly is defined the reduction for left kernel $E^l_X\to (E^l_X)^{red}$.
Note that $E_X^{red},(E^l_X)^{red}$ are naturally embedded into $\cO^{\oplus d}_{X_{red}}(d_{red}-1)$,
where $d_{red}:=\sum d_\al$
\\
\\
These reductions have properties similar to the original kernels, in particular
they fix the \dr s uniquely, in the sense of proposition \ref{Thm.Notions.of.equality of kernels}.
\beR The kernel sheaf of such a reduction can be defined geometrically.
 First, suppose $E_X^{red}$ has a constant rank on the hypersurface.
 Then it defines the rational map:
\beq
X_{red}\stackrel{\phi}{\dashrightarrow}Gr(\P^{rank(E_X^{red})-1},\P^{d-1}),\hspace{1cm} Smooth(X_{red})\ni pt\to Ker (\cM|_{pt})
\eeq
As in the reduced case, let $\tX\norm X_{red}$ be a birational morphism such that $\nu^*\phi$
extends to the morphism on the whole $\tX$. Let $\tau$ be the tautological bundle
on $Gr(\P^{rank(E_X^{red})-1},\P^{d-1})$. Then $\nu^*E_X^{red}(1-d_{red})/Torsion=\phi^*\tau$,
as their fibers coincide. And, moreover, $\phi^*\tau$ determines $E_X^{red}(1-d_{red})$ uniquely.

More generally, let the decomposition $X=\cup X_\al$ satisfy: $X_\al$, $X_\be$ have no common components
for $\al\neq\be$, and
 $E_X^{red}$ has a constant rank on each $X_\al$. Then, as previously, each $E^{red}|_{X_{\al,red}}$ is determined
 uniquely from the corresponding tautological bundle.
\eeR~
\\
\\
The theorem \ref{Thm.Kernel.Sheaves.Properties} is translated almost verbatim:
\bprop\label{Thm.Kernel.Sheaves.Reduced.Locus.Properties}
 Let $X=\cup p_\al X_\al$ be a hypersurface, $deg(X)=d$ and $\cM$ a generically \mg\ \dr.
Let $E_X^{red}$ be the reduction of $E_X$.
 Then
\beq\ber
0\to E_X^{red}\to\cO_{X_{red}}^{\oplus d}(d_{red}-1)\stackrel{\cM}{\to}\cO_{X_{red}}^{\oplus d}(d_{red})\to Coker(\cM)_{X_{red}}^{red}\to0,\\
\eer\eeq
\\{\bf 1.}  The sheaf $E_X^{red}$ is generated by the columns of $\frac{1}{\prod_\al f^{p_\al-1}}\cMv$
as an $\cO_{X_{red}}$ module.
In particular $h^0(E_X^{red})=d$ and $h^0(E_X^{red}(-1))=0$ and
 $h^i(E_X^{red}(j))=0$ for $0<i< n-1$, $j\in\Z$.
\\{\bf 2.}
The sheaves $E_X^{red}$, $Coker(M)_{X_{red}}^{red}$  are torsion free.
\\{\bf 3.} Similar statements hold for $(E^l_X)^{red}$.
\eprop
\bpr  The exactness of the sequences follows immediately from
$\frac{\cMv}{\prod f^{p_\al-1}_\al}\cM=\prod f_\al\one$.
\\{\bf 1.} The proof is identical to that in theorem \ref{Thm.Kernel.Sheaves.Properties}.
\\{\bf 2.} $E_X^{red}$ and $(E^l_X)^{red}$ are torsion free as subsheaves of torsion free sheaves.
For $Coker(M)_{X_{red}}^{red}$, $Coker(M^T)_{X_{red}}^{red}$ the proof is similar to that
in \ref{Thm.Kernel.Sheaves.Properties}.
If $g s=0\in\cO_{X_{red}}$ then $\Big(\frac{1}{\prod_\al f^{p_\al-1}}\cMv\Big) gs=0$ hence $s=\cM s_1$, etc.
\epr
The characterization of kernel sheaves is also translated immediately from the theorem \ref{Thm.Kernel.Sheaves.Classification}.
\bprop\label{Thm.Kernel.Sheaves.Reduced.Locus.Classification}
Let $E_X^{red}$ be a torsion free sheaf on the reduction of the hypersurface $X_{red}=\cup X_\al$, satisfying
$h^0(E_X^{red}(-1))=0$ and $h^i(E_X^{red}(j))=0$ for $0<i<n-1$, $j\in\Z$
and $h^{n-1}(E_X^{red}(1-n))=0$. Suppose the multirank of $E_X^{red}$ is $(p_1,..,p_k)$.
Then $E_X^{red}$ is the restriction of a kernel from $X=\cup p_\al X_\al$,
i.e. $0\to E_X^{red}\to\cO^{\oplus d}_{X_{red}}(d_{red}-1)\stackrel{\cM}{\to}..$
\eprop
\bpr $E_X^{red}$ is a  module over $\cO_{X_{red}}=\cO_{\P^n}/\prod f_\al$.
 Let $\cO_X=\cO_{\P^n}/\prod f^{p_\al}_\al\stackrel{\pi}{\to}\cO_{X_{red}}$ be the natural projection.

Construct an $\cO_X$ module  by defining the action of $\cO_X$ on $E_X^{red}$:
\beq
\rm{for}~ g\in\cO_X ~ \rm{ and} ~ s\in E_X^{red}, ~gs:=\prod_\al f^{p_\al-1}_\al\pi(g)s
\eeq
Denote this module by $E_X$. It is torsion free by construction and all the needed cohomologies vanish.

Hence by theorem \ref{Thm.Kernel.Sheaves.Classification} $E_X$ is the kernel of a \dr, whose reduction is $E_X^{red}$.
\epr

\subsection{Kernels on modifications of curves}\label{Sec.Kernels.on.Modifications}
It is useful to have the classification of kernels in terms of locally free sheaves on the normalization $\tC\norm C$
or, more generally, torsion free sheaves on an intermediate modification $C'\norm C\sset\P^2$.

\subsubsection{When is the pushforward $\nu_*(E)_C$ the kernel sheaf?}
\ \\
Recall (\S\ref{Sec.Preliminaries.XXS.DR.s}) that a torsion free
sheaf $E_C$ is \CCS if $E_C=\nu_*(\nu^*E_C/Torsion)$.
Equivalently, the kernel $E_C$ of $\cM$ is \CCS iff every element of $\cMv$ belongs to the relative adjoint
ideal $Adj_{C'/C}$.

In particular, the cohomology dimensions of these sheaves are preserved under pullbacks and pushfowards.
Hence the theorem \ref{Thm.Kernel.Sheaves.Classification} implies immediate:
\bcor\label{Thm.Kernel.Descent.From.Modification.Classification}
Given a modification $C'\norm C=\cup p_\al C_\al\sset\P^2$ and a torsion free sheaf $E_{C'}$, whose pushforward
$\nu_*(E_{C'})$ is of multi-rank $(p_1,..,p_k)$. The pushforward
$\nu_*(E_{C'})$ is the kernel sheaf of a \dr\ $\cM_C$ iff all the relevant cohomologies vanish:
 $h^0(E_{C'}(-1))=0=h^1(E_{C'}(-1))$.
In this case $\cM_C$ is the unique \CCS \dr\ whose kernel pulls-back to $E_{C'}$.
\ecor
(Here $E_{C'}(-1)=E_{C'}\otimes\nu^*\cO_C(-1)$.)
\bex
1. For the normalization $\tC\norm C$ any torsion free sheaf $E_\tC$ is locally free.
Hence we get the characterization of the locally free sheaves on $\tC$ arising as kernels of \dr s.
\\2. In particular, let $C_d\sset\P^2$ be an irreducible curve, let $\tC_g\norm C_d$ be the normalization.
Then $\nu^*\cO_C(-1)$ is a line bundle of degree $d$. The $\tC/C$-saturated \dr s of $C_d$ correspond
to line bundles of degree $d+g-1$ on $\tC_g$ satisfying: $h^0(\cL_{\tC_g}\otimes\nu^*\cO_C(-1))=0$.
If $g=0$ then the unique such bundle on $\tC=\P^1$ is $\cO_{\P^1}(d-1)$.
Hence there exists unique \dr\ of $C\sset\P^2$ that is $\tC/C$ saturated at each point. More generally,
such line bundles are parameterized by an open dense subset of Jacobian of $C_g$, the complement of
Brill-Noether locus.
\eex
\subsubsection{When is a torsion free sheaf $E_{C'}$ the pullback of some kernel sheaf on $C$?}\
\\
Usually, many torsion free sheaves on $C$ are not pushforwards of locally free sheaves, i.e. are not
of the form $\nu_*E_{C'}$ for any modification $C'\norm C$, cf. \cite{Kerner-Vinnikov2010}.

Thus, while the kernel sheaf has no global sections, $h^0(E_C(-1))=0$, its pull-back can have them.
Let $E_{C'}$ be a torsion free sheaf and $pt\in C'$. Consider the stalk $E_{(C',pt)}$, a module over the local
ring $\cO_{(C',pt)}$. Suppose it is minimally generated by the elements $a_1,..,a_k$. Let
$\nu^{-1}(\cO_{(C,\nu(pt))})\bl a_1,..,a_k\br$ be a vector subspace of $E_{(C',pt)}$.

For generic point $pt\in C'$ one has $\nu^{-1}(\cO_{(C,\nu(pt))})=\cO_{(C',pt)}$, hence
$\nu^{-1}(\cO_{(C,\nu(pt))})\bl a_1,..,a_k\br=E_{(C',pt)}$.
The natural object is the quotient vector space
$E_{(C',pt)}/\nu^{-1}(\cO_{C,\nu(pt)})\bl a_1,..,a_k\br$.
The dimension of this vector space is independent of the choice of generators $\{a_i\}$. We denote this dimension by
 $length_{pt}\quotient{E_{C'}}{\cO_C}$.
\bed
The global section $s\in H^0(E_{C'})$ is said to descend to $C$ locally if for any point $pt\in C'$ there is
a choice of the local generators of $E_{(C',pt)}$ as above,
such that $s\in\nu^{-1}(\cO_{C,\nu(pt)})\bl a_1,..,a_k\br$.
\eed
\bprop\label{Thm.Kernel.Ascent.To.Modification}
The torsion free sheaf $E_{C'}$ is the pull-back of a kernel sheaf on $C$ iff $h^1(E_{C'}(-1))=0$
and $h^0(E_{C'}(-1))\le\sum_{pt\in C'}length_{pt}\quotient{E_{C'}}{\cO_C}$ and no non-zero
global section of $E_{C'}$ descends to $C$ locally.
\eprop
\bpr
$\Rrightarrow$ Suppose $E_C$ is the kernel module, let $E_{C'}=\nu^*(E_C)/Torsion$.
Then $0\to E_{C}\to\nu_*(E_{C'})\to sky\to0$, where $sky$ is
a skyscraper sheaf supported at a finite number of points, at which $E_{C,pt}\subsetneq\nu_*E_{C',pt}$.
Take the cohomology:
\beq
h^0(E_C(-1))=0\to h^0(\nu_*(E_{C'})(-1))\to length(sky)\to h^1(E_C(-1))=0\to h^1(\nu_*(E_{C'})(-1))\to0
\eeq
Hence, $h^1(\nu_*(E_{C'})(-1))=0$ and moreover:
\beq
h^0(\nu_*(E_{C'})(-1))=length(sky)=\sum_{pt\in C}\quotient{\nu_*E_{(C',\nu^{-1}(pt))}}{E_{(C,pt)}}\le\sum_{pt\in C'}\quotient{E_{C'}}{\cO_C}
\eeq
As $h^0(E_C(-1))=0$ no global section of $h^0(\nu_*(E_{C'})(-1))$ descends to $C$.
\\
\\
$\Rrightarrow$ Suppose $E_{C'}$ is given. Take a torsion sheaf, $sky_{C}$, with
$Supp(sky)\sset\cup_{pt\in C}Supp(\nu_*E/\nu^{-1}\cO_{(C,pt)})$ and $length(sky)=h^0(E_{C'}(-1))$
and $length_{pt}(sky)\le length \nu_*E/\nu^{-1}\cO_{(C,pt)}$.
It exists by the assumption $h^0(E(-1))\le\sum_{pt\in C'}length_{pt}\quotient{E_{C'}}{\cO_C}$.

Take a surjection $\nu_*E_{C'}(-1)\to sky\to 0$. Let $E_C$ be the kernel of this map,
\beq
0\to E_{C}\to\nu_*(E_{C'})\to sky\to0
\eeq
Then, by construction, $\chi(E_C(-1))=0$. And by the assumption on non-descent of the global sections of $E_{C'}$:
$h^0(E_C(-1))=0$. Hence $h^1(E_C(-1))=0$.
\epr
\subsection{Families of \dr s}\label{Sec.Kernel.Sheaves.Examples}
Theorem \ref{Thm.Kernel.Sheaves.Properties} and theorem \ref{Thm.Kernel.Sheaves.Classification} translate the
study of \dr s into the study of the sheaves with specific properties. Many related questions are open. We make only a
few remarks.
\li On a smooth curve the generic line bundle of $\deg(\cL(-1))=g-1$ has $h^0(\cL(-1))=0=h^1(\cL(-1))$. Hence the
kernel bundles correspond to an open dense locus on the Jacobian of all the line bundles of the given degree.
\li For a singular, reducible curve $\cup_iC_i$ the $\tC/C$-saturated \dr s correspond to line bundles
$\cL_{\coprod\tC_i}$ with $h^0(\cL_{\coprod\tC_i}\otimes\nu^*\cO_C(-1))=0=h^1(\cL_{\coprod\tC_i}\nu^*\cO_C(-1))$.
 Hence the multi-degree: $\deg(\cL(-1))=\{g(\tC_i)-1\}_i$. Again, the generic line bundle with such a multi-degree
 has no global sections and hence corresponds to the kernel bundle.
\li More generally, the families of \CCS \dr s correspond to some families of torsion free sheaves on $C'$, with
vanishing cohomology.
\li (Continuing example \ref{Ex.Det.Reps.Smooth.Quadrics.P3}.)
A quadric with one $A_1$ singularity is the cone over a smooth plane conic. The kernel bundle corresponds to a line
on the quadric passing through the singular point. The corresponding \dr\ is equivalent to a symmetric one.
\li The \dr s of smooth cubics in $\P^3$ are studied e.g. in \cite{Beauville00}, \cite{Buckley-Košir2007},
for singular cubics cf. \cite{Dolgachev-book}.
\section{Symmetric \dr s}\label{Sec.Symmetric.DetReps}
Here we consider {\em symmetric} local/global \dr s, i.e.
\beq
\cM=\cM^T\in Mat(d\times d,R)\quad \text{ for }R=\cOn\quad \text{ or }R=|\cO_{\P^n}(1)|
\eeq
 In this case the local/global symmetric equivalence is:
$\cM\stackrel{s}{\sim}A\cM A^T$, for $A$- (locally) invertible matrices with entries in $\cOn$ or with constant entries.

As in the ordinary case (property \ref{Thm Localiz Chip off Unity}) the symmetric reduction of \dr s exists.
Choose the local coordinates $(x_1,..,x_n)$ on $(\k^n,0)$.
\bpro
1. Any symmetric matrix $\cM\in Mat(d\times d,\cOn)$ is symmetrically
equivalent to $\one\oplus\cN$, where $\cN=\cN^T$ and $\cN|_{(0,0)}=\zero$. The symmetrically reduced matrix,
i.e. $\cN$, is unique up to the local symmetric equivalence.
\\2. Any local symmetric \dr\  with rational entries is a symmetric reduction of some global symmetric \dr.
If $\cM_1\stackrel{s}{\sim}\cM_2$ locally and $\cM_1$ is the reduction of
$\cM_{global}=\cM^T_{global}$ then $\cM_2$ is also the symmetric reduction of
$\cM_{global}$.
\epro
\bpr
{\bf 1.} (For a part of this statement see also \cite[Lemma 1.7]{Piontkowski2006}.)
Let $\cM\to A\cM A^T$ be a transformation diagonalizing $jet_0\cM$. By further symmetric scaling and permutation of
the rows/columns one
can assume $jet_0(A\cM A^T)=\one\oplus\zero$. Now, in the transformed matrix $A\cM A^T$ kill all the non-constant entries
of the first row/column. This is done symmetrically and decomposes $\cM$ as $\one_{1\times1}\oplus N_1$, where $N_1=N^T_1$.
Do the same for $N_1$ etc.
Uniqueness of the localization is proved as in the property \ref{Thm Localiz Chip off Unity}.
\\
\\{\bf 2.} Let $\cM$ have rational entries, let $g$ be the product of the denominators in all the entries
of $\cM$.  Then $g$ is invertible at the origin and $(g\one)\cM(g\one)$ has polynomial entries.

Given $\cM=\cM^T$ with polynomial entries we should find $\cN=\cN^T$, whose entries are
polynomials of degree at most one, such that $\cN\stackrel{s}{\sim}\one\oplus\cM$.

As in lemma \ref{Thm Localiz Every Local Arises From Global} consider the monomials, appearing in $\cM$,
with the maximal total degree $\deg\cM$. It is enough to show that $\one_{k\times k}\oplus\cM$
is symmetrically equivalent to a matrix $\cN=\cN^T$, such that $\deg\cN\le\deg\cM$ and the number
of monomials in $\cN$, whose total degree is $\deg\cM$, is strictly less than in $\cM$.
Let $x^{a_1}_1..x^{a_n}_n$ be one such monomial, we can assume it is on the diagonal of $\cM$.
Then the step of induction is:
\beq\ber
\bpm 1&0&..&0\\0&1&0&0\\0&0&x^{a_1}_1..x^{a_n}_n+..&..&\\0&0&..&... \epm\rightsquigarrow
\bpm 0&1&..&0\\1&0&0&0\\0&0&x^{a_1}_1..x^{a_n}_n+..&..&\\0&0&..&... \epm\rightsquigarrow
\\\rightsquigarrow\bpm 0&1&x^{a_1}_1..x^{a_i-1}_i..x^{a_n}_n&0\\1&0&0&0\\x^{a_1}_1..x^{a_i-1}_i..x^{a_n}_n&0&x^{a_1}_1..x^{a_n}_n+..&..&\\0&0&..&... \epm
\rightsquigarrow
\bpm 0&1&x^{a_1}_1..x^{a_i-1}_i..x^{a_n}_n&0\\1&0&-\frac{x}{2}&0\\x^{a_1}_1..x^{a_i-1}_i..x^{a_n}_n&-\frac{x}{2}&0+..&..&\\0&0&..&... \epm
\eer\eeq
\epr
A \dr\ is symmetric iff its left and right kernels coincide.
\bpro\label{Thm.Kernel.is.Symmetric.iff.Det.Rep.Is.Symmetric}
Let $\cM_{d\times d}$ be a local/global \dr\ let $E,E^{(l)}\sset R^{\oplus d}$ be its right and
left kernel modules/sheaves. Here $R=\cOn$ in the local case or $R=|\cO_{\P^n}(d-1)|$ in the global case.
\\1. $E\approx E^{(l)}$ iff $\cM$ is equivalent to a symmetric matrix.
\\2. $E=E^{(l)}\sset R^{\oplus d}$ iff $\cM=\cM^T$.
\epro
\bpr The parts $\Lleftarrow$ are obvious.
The converse follows immediately from the uniqueness of minimal free resolution:
\beq
\bM 0&\to&\!\! E&\to& R^{\oplus d}&\stackrel{\cM}{\to}&R^{\oplus d}
\\&&\downarrow \phi&&\downarrow A&&\downarrow B
\\0&\to&E^{(l)}&\to& R^{\oplus d}&\stackrel{\cM^T}{\to}&R^{\oplus d}
\eM
\eeq
An isomorphism $\phi$ induces the isomorphisms $A,B$, by property \ref{Thm.Notions.of.equality of kernels}.
And if $\phi$ is identity then $A,B$ are identities too.
\epr
In the symmetric case the ordinary equivalence implies the symmetric one.
\bprop\label{Thm.Symmetric.Equivalence.iff.Ordinary.Equivalence}
1. If two (local or global) symmetric \dr s of hypersurfaces are (locally or globally) equivalent then they
are symmetrically equivalent.
\\2. Suppose $\cM=\cM^T$ is locally or globally decomposable: $\cM\sim\bpm\cM_1&\zero\\\zero&\cM_2\epm$.
Then there exists a symmetric decomposition: $\cM\stackrel{s}{\sim}\bpm\cN_1&\zero\\\zero&\cN_2\epm$,
 with $\cN^T_\al=\cN_\al\sim\cM_\al$.
\\3. Global symmetric \dr\  decomposes (symmetrically) iff it decomposes locally at the relevant points,
as in Theorem \ref{Thm Decomposability Global from Decomposability Local}.
\eprop
\bpr{\bf Part 1} in the global case is a classical fact e.g. for
two variables see \cite[Chapter VI, \S23,Theorem 3]{Mal'cev-book}. We (re-)prove it both in the
global and local cases.

Let $\cM_1=\cM^T_1$ and $\cM_2=\cM^T_2$ be equivalent, i.e. $\cM_1=A\cM_2B$. We can get rid of $B$
by replacing $\cM_2$ by $(B^T)^{-1}\cM_2B^{-1}$,
so $\cM_2$ stays symmetric. Now one has $\cM_1=A\cM_2=(A\cM_2)^T=\cM_2 A^T$. Note that at each $pt\in X$:
$Ker(A\cM_2|_{pt})=Ker(\cM_2|_{pt})$, thus the comparison of the left and right parts gives:
\beq
\forall pt\in X:\ A^T Ker(\cM_2|_{pt})=Ker(\cM_2|_{pt})\sset\k^d
\eeq
{\em Global case.} Suppose for the generic point of $X$ the vector space $Ker(\cM_2|_{pt})$ is one-dimensional.
Then we get
a constant matrix $A^T$ acting on $\P(\k^d)$, preserving the image $\phi(X)\sset\P(\k^d)$, cf.
\S\ref{Sec.Kernel.Sheaf.Geometric.Definition}. As $Span(\phi(X))$
is the whole ambient space this implies: $A^T$ is diagonal.

In general, for the decomposition $X=\cup p_\al X_\al$,
let $d_\al$ be the generic dimension of $Ker(\cM_2)$ on $X_\al$. Then $A^T\circlearrowright Gr(\P^{d_\al-1},\P^{d-1})$
and $A^T$ preserves $\phi(p_\al X_\al)\sset Gr(\P^{d_\al-1},\P^{d-1})$. Combining this for all the components
 we get again: $A^T$ acts diagonally.

Finally, if $A$ is diagonal we get: $A\cM_2=\cM_2 A$. Thus $A=\tA^2$ with $\cM_1=A\cM_2=\tA\cM_2\tA$.

{\em Local case.} Expand $A$ in powers of local coordinates, then the constant part, $jet_0(A)$, is invertible
and satisfies: $jet_0(A)\cM_2=\cM_2 jet_0(A)^T+higher.order.terms$. Thus, arguing as above we get that $jet_0(A)$
is diagonal and further, is symmetrically equivalent to the identity. So, we assume $jet_0(A)=\one$. Now define
$\sqrt{A}$ as follows. Consider the Taylor series $\sqrt{1+x}:=g(x)$ and define
$\sqrt{A}=\sqrt{\one+(A-\one)}=g(A-\one)$. As $jet_0(A-\one)=\zero$ the series $g(A-\one)$ is well defined
(at least as a formal series).
From $A\cM_2=\cM_2A^T$ we get $A^j\cM_2=\cM_2(A^T)^j$, thus $\sqrt{A}\cM_2=\cM_2\sqrt{A^T}$. Therefore we get:
\beq
\cM_1=A\cM_2=\sqrt{A}\cM_2\sqrt{A^T}
\eeq

{\bf Part 2.} As $E=E^{(l)}$ their restrictions to the components $(X_\al,0)$ coincide too.
Hence $\cM\sim\cN_1\oplus\cN_2$ where $\cN^T_\al=\cN_\al$.
Now, by the first statement we get that the equivalence can be chosen symmetric.
\\{\bf Part 3} follows from the ordinary decomposability and part 2.
\epr
Finally we characterize  the sheaves that are kernels of symmetric \dr s.
\bprop\label{Thm Symmetric Case.Classification.Kernel Sheaves}
Let $E_X$ be a torsion-free sheaf on the hypersurface $X\sset\P^n$.
Assume all the relevant cohomologies vanish (as in theorem \ref{Thm.Kernel.Sheaves.Classification}) and $E_X\approx E_X^*\otimes w_X(n)=Hom(E_X,w_X(n))$.
Then there exists $\cM^T=\cM\in Mat(d\times d,H^0(\cO_{\P^n}(1)))$ such that
\beq
0\to E_X\to\cO_X^{\oplus d}(d-1)\stackrel{\cM}{\to}\cO_X^{\oplus d}(d)\to..
\eeq
If moreover $X$ is reduced and $E_X$ is locally free then $n\le2$.
\eprop
\bpr The main construction is done in the proof of theorem \ref{Thm.Kernel.Sheaves.Classification},
it realizes $E_X$ and $E^l_X=E_X^*\otimes w_X(n)$ as the right and left kernels of $\cM$.
Hence, as $E_X\approx E^l_X$ we get from proposition
\ref{Thm.Kernel.is.Symmetric.iff.Det.Rep.Is.Symmetric} that $\cM$ is equivalent to a symmetric matrix.

The last statement is proved as the last statement of theorem \ref{Thm.Kernel.Sheaves.Properties}.
Local freeness implies that corank of $\cM_X$ is one at any point. But in the parameter space of all the $d\times d$
symmetric matrices the subset of matrices of corank at least two is of codimension three. Hence if $n\ge3$ there
will be always a point $pt\in X$ such that $dim(Ker\cM|_{pt})>1$.
\epr
\beR
1. Note that even if the hypersurface is reduced and the kernel is of rank one and generically locally free,
we ask for $E_X\approx E_X^*\otimes w_X(n)$ rather than $E_X\otimes E_X\approx w_X(n)$
as $E_X\otimes E_X$ is not a torsion free sheaf, unless $E_X$ is locally free. In fact, even with the torsion factored out,
the map $E_X\otimes E_X/Torsion\to w_X(n)$ is injective (as its kernel would be a torsion subsheaf) but
is not surjective.
\\2. In the case of curves the "self-dual" sheaves of the proposition are the well known theta characteristics,
 see \cite{Harris-1982} for the smooth case and \cite{Piontkowski2007} for reduced case. In particular,
 for any collection of the "local types", i.e. the stalks of $E_X$ at the points where $E_X$ is not locally free,
 there exists a theta characteristic with this collection of types.
\eeR

\section{Self-adjoint \dr s}\label{Sec.Self.Adjoint.dr s}
In this section we work over $\C$, i.e. $\P^n=\P^n_\C$. Let $\P^n_\R$ be the real projective space .
For a subscheme $X\sset\P^n$ we denote the set of its real points by $X_\R$.

\subsection{Setup}
Let the hypersurface $X\sset\P^n$ be defined over $\R$, i.e. its defining polynomial has real coefficients.
Let $\tau\circlearrowright\C$ be the complex conjugation, so $\tau$ acts on $\P^n$ and on $X$.
Thus $\tau$ acts on the twisting sheaves $\tau\circlearrowright\cO_{\P^n}(d)$ and
$\tau\circlearrowright\cO_X(d)$. The complex conjugation acts on the set of all the sheaves of embedded modules.
\bed
Let $E_X\sset \oplus_\al\cO_X(d_\al)$ be a sheaf of modules. The induced conjugation,
$E_X^\tau:=\tau(E_X)$, is defined by the action on sections, i.e. each section is sent to its conjugate.
\eed
For any matrix define $A^\tau:=\bar{A}^T$, i.e. both conjugated and transposed.
\bed
A local or global \dr\  is called self-adjoint if for any $pt\in\P^n_\R$ one has: $\cM|_{pt}=\cM|_{pt}^\tau$.
Self-adjoint \dr s are considered up to Hermitean equivalence: $\cM\stackrel{\tau}{\sim}A\cM A^\tau$.
\eed
As in the symmetric case, the self-adjointness can be expressed in terms of kernels and the
ordinary equivalence implies an almost Hermitean equivalence,
generalizing \cite[\S9, theorem 8]{Vinnikov89}:
\bprop\label{Thm.Self-adjoint.Equivalence.Iff.Ordinary.Equiv}
1. $\cM=\cM^\tau$ iff $E^\tau=E^{(l)}\sset\cO_X^{\oplus d}(d-1)$. And $\cM\sim\cM^\tau$ iff $E^\tau\approx E^{(l)}$.
\\2. Let $\cM$, $\cM'$ be (local or global) \dr s of the same hypersurface.
If $\cM\sim\cM'$ and both are self-adjoint then either $\cM\stackrel{\tau}{\sim}\cM'$
or $\cM\stackrel{\tau}{\sim}-\cM'$ or
\beq
\cM\stackrel{\tau}{\sim}\bpm \cM_1&\zero\\\zero&\cM_2\epm \text{ and }
\cM'\stackrel{\tau}{\sim}\bpm \cM_1&\zero\\\zero&-\cM_2\epm
\eeq
\eprop
\bpr Part one is obvious.
\\Part two. As in proposition \ref{Thm.Symmetric.Equivalence.iff.Ordinary.Equivalence},
starting from $\cM=A\cM'=(A\cM')^\tau=\cM' A^\tau$ we get: $A$ preserves the embedded kernel at each point of $X$.
If $A$ is constant (e.g. in the global case) then $A$ is diagonal, so after reshuffling the rows/columns and
rescaling (over reals!) we can assume $A=\bpm\one&\zero\\\zero&-\one\epm$. This implies that $\cM'$ is block diagonal
and gives the statement in the global case.

In the local case we get $jet_0(A)=\bpm\one&\zero\\\zero&-\one\epm$, hence the expansion of $jet_0(A)\times A$
in local coordinates begins with the constant term $\one$. Therefore we can use the Taylor expansion
of $\sqrt{jet_0(A)\times A}=\sqrt{\one+(jet_0(A)\times A-\one)}$ and proceed as in the proof of
proposition \ref{Thm.Symmetric.Equivalence.iff.Ordinary.Equivalence}.
\epr

\subsection{Classification}
The classification of self-adjoint \dr s of smooth plane curves was done in \cite[\S9, theorem 7]{Vinnikov89}.
As in the ordinary/symmetric cases we generalize to an arbitrary hypersurface.
\bthe Let $X\sset\P^n$ be an arbitrary hypersurface, defined over $\R$.
A torsion free sheaf $E_X$ is the kernel sheaf of a self-adjoint \dr\ iff:
all the relevant cohomologies vanish (as in theorem \ref{Thm.Kernel.Sheaves.Classification})
and $E_X^{(l)}=E_X^*\otimes w_X(n)\approx E_X^\tau$.
\ethe
\bpr The statement on vanishing cohomologies was proved in theorems
\ref{Thm.Kernel.Sheaves.Properties} and \ref{Thm.Kernel.Sheaves.Classification}.

For the self-conjugacy of the kernel sheaf, the direction $\Rrightarrow$ is obvious.
The inverse direction is proved exactly as in proposition \ref{Thm Symmetric Case.Classification.Kernel Sheaves}.
We should only check Step 3 of the original construction.

Let $E_X\isom(E_X^{(l)})^\tau$ be the isomorphism of sheaves.
By property \ref{Thm.Notions.of.equality of kernels} it extends to the global automorphism
 $\phi\circlearrowright\cO_X^{\oplus d}(d-1)$. As $\phi$ is a global automorphism it is
 presented by a constant invertible matrix. Choose the basis of $\cO_X^{\oplus d}(d-1)$ so
that $\phi$ becomes identity. Namely $E_X$ and $(E_X^{(l)})^\tau$ coincide as embedded sheaves.
Finally, once the sheaves are identified, choose the same bases: $E_X,(s_1,..,s_d)=(E_X^{(l)})^\tau,(s_1,..,s_d)$.
\epr
\subsection{An application to hyperbolic polynomials}\label{Sec.Self Adjoint dr of hyperbolic polynomials}
Let $\cM\in Mat\Big(d\times d,H^0(\cO_{\P^n_\R}(1))\Big)$ be a self-adjoint positive definite \dr\ of a real
 projective hypersurface.
Namely, $\cM$ is a matrix of linear forms in homogeneous coordinates on $\P^n_\R$ and for any $pt\in\P^n_\R$
one has: $\cM|_{pt}=\cM|_{pt}^\tau$ and at least for one point
$\cM|_{pt}$ is positive definite.
\bprop
Let $\cM$ be positive definite at least at one point. Then for some choice of coordinates on $\P^n_\R$:
$\cM=\sum^n_{i=0} \cM_i x_i$, where all the matrices $\cM_i$ are positive definite.
\eprop
\bpr If $\cM$ is positive definite at $pt\in\P^n_\R$ then it is positive definite on some open neighborhood
of this point.
Let $pt_0,..,pt_n\in\P^n_\R$ be some points close to $pt$, such that $\cM$ is positive definite at these
points and the points do not lie in a hyperplane, i.e. they span $\R^n$ locally.

Consider $\P^n_\R$ as $Proj(\R^{n+1})$, so one can choose the coordinate axes $\hat{x}_0$, .., $\hat{x}_n$
of $\R^{n+1}$, corresponding to these points. Hence $\cM|_{pt_j}= x_j\cM_j$ for $x_j>0$. Hence
in this coordinate system all $\cM_i$ are positive definite.
\epr
\bed
A homogeneous polynomial $f\in\R[x_0,..,x_n]$ is called hyperbolic
if there exists a point $pt\in\P^n_\R$ with the property: for any line $L$ through $pt$
all the complex roots of $f|_L=0$ are real. (In other words the line intersects the
hypersurface $\{f=0\}\sset\P^n_\R$ at $\deg(f)$ real points, counted with multiplicities.)
\eed
For the general introduction to the theory of hyperbolic polynomials cf. \cite{Gårding1951}, \cite{Gårding1959},
\cite{Lax1958}, \cite{Guler1997}, \cite{Bauschke-Guler-Lewis-Sendov2001}, \cite{Renegar2006}.
For a given hyperbolic polynomial the union of all points satisfying the definition above is
called {\em the region of hyperbolicity}. It is a convex set in $\P^n_\R$ under projection $\R^{n+1}\to\P^n_\R$,
 i.e. its preimage in $\R^{n+1}$
is the disjoint union of two convex sets.
Suppose the real homogeneous polynomial in three variables defines a smooth plane curve. Then the polynomial
is hyperbolic iff the corresponding curve has the maximal possible number of nested ovals.

A self-adjoint positive-definite \dr\ defines a hyperbolic polynomial
(see \cite[\S6]{Vinnikov93} for curves and \cite{Helton-Vinnikov2007} for hypersurfaces).
In this case the hyperbolicity region consists of all points $pt\in\P^n_\R$ such that $\cM|_{pt}$ is (positive
or negative) semi-definite (see the citations above).

A naive converse statement could be: if $f(x_0,..,x_n)$ is a hyperbolic polynomial with the hyperbolicity region $\R^n_{\ge0}$
 then $f^N$ has a positive-definite \dr\ for $N\gg0$.  This turns to be wrong in higher dimensions, \cite{Brändén2010}.
 But for $n=2$ this is true even for $f$ itself, not only for its higher multiples, \cite{Dubrovin1983}, \cite{Vinnikov93},
 \cite{Helton-Vinnikov2007}.
  A weaker statement, {\em there exists a hyperbolic polynomial $g(x_0,..,x_n)$ with the hyperbolicity region
  containing that of $f$, such that $gf$ is determinantal}, has not been checked yet.

\bthe\label{Thm.Hyperbol.Curve.Has.Smooth.Branches.Det.Rep.is.Maximal}
1. Let $X=\cup p_\al X_\al\sset\P^n$ be a hypersurface defined over $\R$ by a hyperbolic polynomial.
The real part of its reduced locus, $\cup X^\R_\al$, can have at most one (real) singular point with a non-smooth
locally irreducible component. In the later case the region of hyperbolicity degenerates to this singular point.
\\2. In particular, if the hypersurface $X$ possesses a self-adjoint positive definite \dr\
then the germ of its reduced locus at each of
its (real) singular point is the union of smooth hypersurfaces. Hence there exists a finite modification
$\tX=\coprod p_\al\tX_\al\to X=\cup X_\al$ such that the reduced (real) locus $\coprod\tX^\R_\al$ is smooth.
\\3. Let $\cM$ be a self-adjoint positive-definite representation of the hypersurface $X\sset\P^n$.
Let $\tX\to X$ be the finite modification as above.
 The representation is $\tX/X$-saturated at all the real points of $X$.
\ethe
For the definition of $\tX/X$-saturated see \S\ref{Sec.Preliminaries.XXS.DR.s}.
\bpr {\bf 1.} Let $pt\in\P^n_\R$ be a hyperbolic point. Suppose $X_\R$ has a singular locally irreducible component
 $(X^\R_\al,0)\sset(X_\R,0)$, i.e. the multiplicity $mult(X^\R_\al,0)>1$. Suppose $pt\neq 0$.
 We prove that there exists a family of real lines $L_1(t)$, all passing through the point $pt$ and satisfying:
\ls $L_1(0)\ni 0$, hence $\deg(L_1(0)\cap (X^\R_\al,0))\ge mult(X^\R_\al,0)$.
\ls There exists a small neighborhood of $0\in\P^n_\R$, in the classical topology, such that for any $t\neq0$ the total intersection
degree $\deg(L_1(t)\cap X^\R_\al)$ in this neighborhood is less than $mult(X^\R_\al,0)$.

This will contradict hyperbolicity of the polynomial, implying that either $0=pt$ or $(X^\R_\al,0)$ is smooth.

First, the problem can be reduced to the planar case, i.e. $n=2$. Indeed, let $L_2\sset(\R^n,0)$ be
the generic real two-dimensional plane through the origin. Then $(L_2\cap X^\R_\al,0)$ is locally irreducible
and singular. Hence, it is enough to present a family of lines in $L_2$, with the needed properties relatively to
the curve $(L_2\cap X^\R_\al,0)$.

Now, observe that a locally irreducible real plane curve $(L_2\cap X^\R_\al,0)$ divides the small neighborhood
(in classical topology) of $0\in L_2$ into {\em two} parts. Hence, if $mult(L_2\cap X^\R_\al,0)=mult(X^\R_\al,0)>2$
then any family of lines $L_1(t)\sset L_2$, whose generic fiber does not pass through the origin and $0\in L_1(0)$,
has the needed properties.

If $mult(L_2\cap X^\R_\al,0)=mult(X^\R_\al,0)=2$ then the singularity is of type $A_k$, i.e. after a (real-analytic)
change of coordinates the curve is defined by $y^2=\pm x^{k+1}$. Now the family is constructed directly.
\\
\\{\bf 2.} If the hypersurface $X\sset\P^n$ has a self adjoint positive definite \dr\ (and not just semi-definite)
then the
hyperbolic region is a non-empty open set. Hence by the previous part the hypersurface
 has no singular points with singular
locally irreducible components. Hence, if $X=\cup p_\al X_\al$, there exists the finite modification
$\coprod p_\al\tX_\al\to\cup p_\al X_\al$, where each $\tX^\R_\al\to X^\R_\al$ is the normalization
and $\tX^\R_\al$ is smooth.
\\
\\{\bf 3.}
As the hypersurface has only smooth multiple local components, it is enough to prove that the representation
is locally completely decomposable at every real singular point of $X$, according to the decomposition
into the distinct multiple components, see property \ref{Thm.X'/X-sat.are.decomposable}.

Let $0\in X_\R$ be a real singular point. Let $pt\in\P^n_\R\smin X_\R$ be
a point in the hyperbolic region. We can assume the line $\overline{(pt,0)}$ is not tangent to $X_\R$ at any point.
Let $\ep\in\P^n_\R$ be a point near $0\in X_\R$.
Consider the line $\overline{(pt,\ep)}$. We can assume this line does not lie on $X_\R$.

Restrict the original \dr\ to this line. So, we get a one-dimensional family of  matrices:
\beq
(1-t)\cM|_\ep+t\cM|_{pt}
\eeq
By construction $\det\Big((1-t)\cM|_\ep+t\cM|_{pt}\Big)$ vanishes
precisely for those values of $t$ where the line intersects $X_\R$. By hyperbolicity there
are $deg(X)$ such points (counted with multiplicities).

As $\cM|_{pt}$ is self-adjoint and positive definite, it can be presented as $\cM|_{pt}=U_{pt} U^\tau_{pt}$.
As $pt\in\P^n_\R\smin X$ one gets that $U_{pt}$ is invertible.
Hence the equation above can be presented as:
\beq
\det\Big(t\one+(1-t)U_{pt}^{-1}\cM|_\ep U_{pt}^{-\tau}\Big)=0
\eeq
i.e. as the equation for the eigenvalues of a matrix. As the matrix is self-adjoint, the corresponding
eigenvectors are orthogonal. As $\ep\to0$ they converge to the mutually orthogonal
eigenvectors of $U_{pt}^{-1}\cM|_0 U_{pt}^{-\tau}$, see e.g. \cite{Baumgärtel}.

Hence we obtain: in the limit $\ep\to0$ the normalized sections of the kernel of $\cM|_\ep$ have
 linearly independent limits. Hence, by property \ref{Thm.Kernel.Has.Independent.Fibres.Decomposable.},
the \dr\ $\cM$ is completely decomposable near $0\in X_\R$. Thus it is $\tX/X$-saturated.
\epr

\end{document}